\documentclass[10pt]{article} 
\usepackage{fancyhdr}
\usepackage{amsmath, amssymb, latexsym, amscd, amsthm,amsfonts,amstext}
\usepackage{subfigure}
\usepackage{enumerate}
\usepackage{xcolor}
\usepackage{textcomp}
\usepackage{mathtools}

\def\div{\operatorname{div}}
\usepackage{graphicx}
\newtheorem{theorem}{Theorem}[section]
\newtheorem{lemma}[theorem]{Lemma}
\newtheorem{proposition}[theorem]{Proposition}
\newtheorem{corollary}[theorem]{Corollary}
\numberwithin{equation}{section}
\newtheorem{definition}[theorem]{Definition}
\newtheorem{remark}[theorem]{Remark}

\newenvironment{indeed}[0]{Indeed, }{}

\numberwithin{equation}{section}

\usepackage[mathscr]{eucal}
\usepackage{graphicx}
\usepackage{subfigure}
\usepackage{fancyhdr}
\usepackage{enumerate}
\usepackage{xcolor}
\def\diam{\operatorname{diam}}

\usepackage{tikz}

\usepackage{url}

\def\div{\operatorname{div}}
\allowdisplaybreaks

\title{The equivalent medium for the elastic scattering 
by many small rigid bodies and applications} 

\author{
{Fadhel  Al-Musallam}
\thanks{Department of Mathematics, 
Kuwait University, P.O. Box 13060, Safat, Kuwait.}
\and 
{Durga Prasad Challa$^{\ddag}$}
\thanks{Department of Mathematics, Tallinn University of Technology,  Tallinn, Estonia. {durga.challa@ricam.oeaw.ac.at}}
\and
{Mourad Sini}
\thanks{$^\dag$ Radon Institute (RICAM), Austrian Academy of Sciences, Altenbergerstrasse 69, A4040, Linz, Austria.}
}
\begin{document}
\maketitle

\begin{abstract}
{ We deal with the elastic scattering by a large number $M$ of rigid bodies, $D_m:=\epsilon B_m+z_m$, of arbitrary shapes with  $ 0<\textcolor{black}{\epsilon}<<1$ and with constant Lam\'e coefficients $\lambda$ and $\mu$. 
 We show that, when these rigid bodies are distributed arbitrarily (not necessarily periodically) in a bounded region $\Omega$ of $\mathbb{R}^3$
 where their number is $M:=M(\textcolor{black}{\epsilon}):=O(\textcolor{black}{\epsilon}^{-1})$ and the minimum distance between them is $d:=d(\textcolor{black}{\epsilon})\approx \textcolor{black}{\epsilon}^{t}$ with $t$ in some appropriate range,
 as $\textcolor{black}{\epsilon} \rightarrow 0$, the generated far-field patterns approximate the far-field patterns generated by an equivalent medium given by $\omega^2\rho I_3-(K+1)\mathbf{C}_0 $ where $\rho$ is the density of the background medium (with $I_3$ as the unit matrix) and $(K+1)\mathbf{C}_0$ is the shifting (and possibly variable) coefficient. 
 This shifting coefficient is described by the two coefficients $K$ and $\mathbf{C}_0$ (which have supports in $\overline{\Omega}$) modeling the local distribution of the small bodies and their geometries, respectively.   
 In particular, if the distributed bodies have a uniform spherical shape then the equivalent medium is isotropic while for general shapes it might be anisotropic (i.e. $\mathbf{C}_0$ might be a matrix).
  In addition, if the background density $\rho$ is variable in $\Omega$ and $\rho =1$ in $\mathbb{R}^3\setminus{\overline{\Omega}}$, then if we remove from $\Omega$   appropriately distributed small bodies then the equivalent medium will be equal to $\omega^2 I_3$ in $\mathbb{R}^3$, i.e. the obstacle $\Omega$ characterized by $\rho$ is approximately cloaked at the given and fixed frequency $\omega$.    
}
{Elastic wave scattering, Small-scatterers, Effective medium.}
\\
2000 Math Subject Classification: 35J08, 35Q61, 45Q05
\end{abstract}

\section{Introduction and statement of the results}\label{Introduction-smallac-sdlp}
\subsection{The background}

 Let $B_1, B_2,\dots, B_M$ be $M$ open, bounded and simply connected sets in $\mathbb{R}^3$ with Lipschitz boundaries,
containing the origin.
We assume that their sizes and Lipschitz constants are uniformly bounded.
We set $D_m:=\epsilon B_m+z_m$ to be the small bodies characterized by the parameter
$\epsilon>0$ and the locations $z_m\in \mathbb{R}^3$, $m=1,\dots,M$.
\par
 Assume that the Lam\'e coefficients $\lambda$ and $\mu$ are constants satisfying $ \mu > 0 \mbox{ and } 3\lambda+2\mu >0$ and the mass density $\rho$ to be a constant that we normalize to a unity.
Let $U^{i}$ be a solution of the Navier equation $(\Delta^e + \omega^{2})U^{i}=0 \mbox{ in } \mathbb{R}^{3}$, $\Delta^{e}:=(\mu\Delta+(\lambda+\mu)\nabla \div)$.
We denote by  $U^{s}$ the elastic field scattered by the $M$ small bodies $D_m\subset \mathbb{R}^{3}$ due to
the incident field $U^{i}$. We restrict ourselves to the scattering by rigid bodies. Hence the total field $U^{t}:=U^{i}+U^{s}$
satisfies the following exterior Dirichlet problem of the elastic waves
\begin{equation}
(\Delta^e + \omega^{2})U^{t}=0 \mbox{ in }\mathbb{R}^{3}\backslash \left(\mathop{\cup}_{m=1}^M \bar{D}_m\right),\label{elaimpoenetrable}
\end{equation}
\begin{equation}
U^{t}|_{\partial D_m}=0,\, 1\leq m \leq M\label{elagoverningsupport}
\end{equation}

with the Kupradze radiation conditions (K.R.C)
\begin{equation}\label{radiationcela}
\lim_{|x|\rightarrow\infty}|x|(\frac{\partial U_{p}}{\partial|x|}-i \textcolor{black}{\kappa_{p}}U_{p})=0,  \mbox{  and  }
\lim_{|x|\rightarrow\infty}|x|(\frac{\partial U_{s}}{\partial|x|}-i \textcolor{black}{\kappa_{s}}U_{s})=0,
\end{equation}
where the two limits are uniform in all the directions $\hat{x}:=\frac{x}{|x|}\in \mathbb{S}^{2}$.
Also, we denote $U_{p}:=-\textcolor{black}{\kappa_{p}^{-2}}\nabla (\nabla\cdot U^{s})$ to be the longitudinal
(or the pressure or P) part of the field $U^{s}$ and $U_{s}:=\textcolor{black}{\kappa_{s}^{-2}}\nabla\times(\nabla\times U^{s})$ to be the transversal (or the shear or S) part of
the field $U^{s}$ corresponding to the Helmholtz decomposition $U^{s}=U_{p}+U_{s}$. The constants $\textcolor{black}{\kappa_{p}}:=\frac{\omega}{c_p}$ and
$\textcolor{black}{\kappa_{s}}:=\frac{\omega}{c_s}$ are known as the longitudinal and transversal wavenumbers, $c_p:=\sqrt{\lambda+2\mu}$ and $c_s:=\sqrt{\mu}$ are the corresponding phase velocities, respectively and $\omega$ is the frequency.
\par

The scattering problem (\ref{elaimpoenetrable}-\ref{radiationcela}) is well posed in the H\"{o}lder or Sobolev spaces,
see \cite{C-K:1998,C-K:1983,Kupradze:1965, K-G-B-B:1979} for instance, and the scattered field $U^s$ has the following asymptotic expansion:

\begin{equation}\label{Lamesystemtotalfieldasymptoticsmall}
 U^s(x) := \frac{e^{i\textcolor{black}{\kappa_{p}}|x|}}{|x|}U^{\infty}_{p}(\hat{x}) +
\frac{e^{i\textcolor{black}{\kappa_{s}}|x|}}{|x|}U^{\infty}_{s}(\hat{x}) + O(\frac{1}{|x|^{2}}),~|x|\rightarrow \infty
\end{equation}
uniformly in all directions $\hat{x}\in \mathbb{S}^{2}$. The longitudinal part of the far-field, i.e. $U^{\infty}_{p}(\hat{x})$ is normal to $\mathbb{S}^{2}$
while the transversal part $U^{\infty}_{s}(\hat{x})$ is tangential to $\mathbb{S}^{2}$. We set $U^{\infty}:=(U_p^{\infty}, U_s^{\infty})$. 
 
As usual, we use plane incident waves of the form
$U^{i}(x,\theta):=\alpha\theta\,e^{i\textcolor{black}{\kappa_{p}}\theta\cdot x}+\beta\theta^\bot\,e^{i\textcolor{black}{\kappa_{s}}\theta\cdot x}$,
where $\theta^{\bot}$ is any direction in $\mathbb{S}^{2}$ perpendicular to the incident direction $\theta \in \mathbb{S}^2$, $\alpha,\beta$
are arbitrary constants. The functions
$U^{\infty}_p(\hat{x}, \theta):=U^{\infty}_p(\hat{x})$ and $U^{\infty}_s(\hat{x}, \theta):=U^{\infty}_s(\hat{x})$ for $(\hat{x}, \theta)\in \mathbb{S}^{2} \times\mathbb{S}^{2}$
are called the P-part and the S-part of the far-field pattern respectively.
\newline

\begin{definition} 
\label{Def1}
We define
\begin{enumerate}

 \item  $d:=\min\limits_{\substack{m\neq j\\1\leq m,j\leq M }} d_{mj},$
$\text{where}\,d_{mj}:=dist(D_m, D_j)$. We assume that $
0\,<\,d\,\leq\,d_{\max},$
and $d_{\max}$ is given.
\item $\omega_{\max}$ as the upper bound of the used frequencies, i.e. $\omega\in[0,\,\omega_{\max}]$.
\item $\Omega$ to be a bounded domain in $\mathbb{R}^3$ containing the small bodies $D_m,\,m=1,\dots,M$.
\end{enumerate}
\end{definition}
\bigskip

\subsection{The results for a homogeneous elastic background}
We assume that $D_m=\epsilon{B}_m+z_m, m=1,\dots,M$, with the maximal diameter $a:=\epsilon \max_{m} \diam (B_m)$, are non-flat Lipschitz obstacles, i.e. $D_m$'s are Lipschitz obstacles and there exist constants $t_0 \in (0, 1]$ such that
$
 B^{3}_{t_0\frac{a}{2}}(z_m)\subset\,D_m\subset\,B^{3}_{\frac{a}{2}}(z_m),
$ where $t_m $ are assumed to be uniformly bounded from below by a positive constant. In \cite{C-S:2015}, we have shown that there exist two positive constants $\textcolor{black}{\epsilon}_0$ and $c_0$ depending only on the size of $\Omega$, the
Lipschitz character of $B_m,m=1,\dots,M$, $d_{\max}$ and $\omega_{\max}$ such that
if
\begin{equation}\label{conditions-elasma}
\textcolor{black}{\epsilon} \leq \textcolor{black}{\epsilon}_0 ~~ \mbox{and} ~~ \sqrt{M-1}\frac{\textcolor{black}{\epsilon}}{d}\leq c_0
\end{equation}
 then we have the following asymptotic expansion
 for the P-part, $U^\infty_p(\hat{x},\theta)$, and the S-part, $U^\infty_s(\hat{x},\theta)$, of the far-field pattern:
 \begin{eqnarray}
  &\hspace{-1cm}U^\infty_p(\hat{x},\theta)=&\frac{1}{4\pi\,c_p^{2}}(\hat{x}\otimes\hat{x})\left[\sum_{m=1}^{M}e^{-i\frac{\omega}{c_p}\hat{x}\cdot z_{m}}Q_m\right.\nonumber\\
  &&\hspace{-2cm}\left.+O\left(M\left[{\textcolor{black}{\epsilon}}^2+\frac{{\textcolor{black}{\epsilon}}^3}{d^{5-3\alpha}}+\frac{ {\textcolor{black}{\epsilon}}^4}{d^{9-6\alpha}}\right]+M(M-1)\left[\frac{{\textcolor{black}{\epsilon}}^3}{d^{2\alpha}}+\frac{{\textcolor{black}{\epsilon}}^4}{d^{4-\alpha}}+
  \frac{{\textcolor{black}{\epsilon}}^4}{d^{5-2\alpha}}\right]+M(M-1)^2\frac{{\textcolor{black}{\epsilon}}^4}{d^{3\alpha}}\right) \right],
   \label{x oustdie1 D_m farmainp-near}\\
 &\hspace{-1cm}U^\infty_s(\hat{x},\theta)=& \frac{1}{4\pi\,c_s^{2}}({\rm \textbf{I}}- \hat{x}\otimes\hat{x})\left[\sum_{m=1}^{M}e^{-i\frac{\omega}{c_s}\hat{x}\cdot\,z_m}Q_m\right.\nonumber\\
 &&\hspace{-2cm}\left.+O\left(M\left[{\textcolor{black}{\epsilon}}^2+\frac{{\textcolor{black}{\epsilon}}^3}{d^{5-3\alpha}}+\frac{ {\textcolor{black}{\epsilon}}^4}{d^{9-6\alpha}}\right]+M(M-1)\left[\frac{{\textcolor{black}{\epsilon}}^3}{d^{2\alpha}}+\frac{{\textcolor{black}{\epsilon}}^4}{d^{4-\alpha}}+
 \frac{{\textcolor{black}{\epsilon}}^4}{d^{5-2\alpha}}\right]+M(M-1)^2\frac{{\textcolor{black}{\epsilon}}^4}{d^{3\alpha}}\right) \right],
 \label{x oustdie1 D_m farmains-near}
  \end{eqnarray}
 where $\alpha$, $0<\alpha\leq1$, is a parameter describing the relative distribution of the small bodies.

 The vector coefficients $Q_m$, $m=1,..., M,$ are the solutions of the following linear algebraic system
\begin{eqnarray}\label{fracqcfracmain}
 C_m^{-1}Q_m &=&-U^{i}(z_m, \theta)-\sum_{\substack{j=1 \\ j\neq m}}^{M} \Gamma^{\omega}(z_m,z_j)Q_j,~~
\end{eqnarray}
for $ m=1,..., M,$ with $\Gamma^{\omega}$ denoting the Kupradze matrix of the fundamental solution to the Navier equation with frequency $\omega$, $C_m:=\int_{\partial D_m}\sigma_m(s)ds$ and $\sigma_{m}$ is
the solution matrix of the integral equation of the first kind
\begin{eqnarray}\label{barqcimsurfacefrm1main}
\int_{\partial D_m}\Gamma^{0}(s_m,s)\sigma_{m} (s)ds&=&\rm \textbf{I},~ s_m\in \partial D_m,
\end{eqnarray}

with $\rm \textbf{I}$ the identity matrix of order 3.
 
\bigskip

Consider now the special case $d_{min} {\textcolor{black}{\epsilon}}^t \leq d\leq d_{max}{\textcolor{black}{\epsilon}}^t$ and $M\leq M_{max} {\textcolor{black}{\epsilon}}^{-s}$ with  $t,s>0$, $d_{min}$, $d_{max}$ and $M_{max}$ are positive. 
Then the asymptotic expansions (\ref{x oustdie1 D_m farmainp-near}-\ref{x oustdie1 D_m farmains-near}) can be rewritten as
\begin{eqnarray}
 U^\infty_p(\hat{x},\theta)&=&\frac{1}{4\pi\,c_p^{2}}(\hat{x}\otimes\hat{x})\big[\sum_{m=1}^{M}e^{-i\frac{\omega}{c_p}\hat{x}\cdot z_{m}}Q_m \nonumber \\
 &&+O\left({\textcolor{black}{\epsilon}}^{2-s}+{\textcolor{black}{\epsilon}}^{3-s-5t+3t\alpha}+{\textcolor{black}{\epsilon}}^{4-s-9t+6t\alpha}+{\textcolor{black}{\epsilon}}^{3-2s-2t\alpha}+{\textcolor{black}{\epsilon}}^{4-3s-3t\alpha}+{\textcolor{black}{\epsilon}}^{4-2s-5t+2t\alpha}\right)\big],\quad \label{x oustdie1 D_m farmainp-near*}\\
 U^\infty_s(\hat{x},\theta)&=& \frac{1}{4\pi\,c_s^{2}}({\rm \textbf{I}}- \hat{x}\otimes\hat{x})\big[\sum_{m=1}^{M}e^{-i\frac{\omega}{c_s}\hat{x}\cdot\,z_m}Q_m \nonumber \\
 &&+O\left({\textcolor{black}{\epsilon}}^{2-s}+{\textcolor{black}{\epsilon}}^{3-s-5t+3t\alpha}+{\textcolor{black}{\epsilon}}^{4-s-9t+6t\alpha}+{\textcolor{black}{\epsilon}}^{3-2s-2t\alpha}+{\textcolor{black}{\epsilon}}^{4-3s-3t\alpha}+{\textcolor{black}{\epsilon}}^{4-2s-5t+2t\alpha}\right) \big].\quad \label{x oustdie1 D_m farmains-near*}
 \end{eqnarray}
 
As $\textcolor{black}{\epsilon}\rightarrow0$, the error term tends to zero for $t$ and $s$ such that 
\begin{equation}\label{general-condition-s-t}
0<t<1 \mbox{ and } 0<s<\min\{2(1-t),\,3-5t+3t\alpha,2-\frac{5}{2}t+t\alpha,4-9t+6t\alpha, \frac{3}{2}-t\alpha,\frac{4}{3}-t\alpha\}.
\end{equation}

 In \cite{C-S:2015}, we have shown that $Q_m\approx \textcolor{black}{\epsilon}$, then we have 
 the upper bound
\begin{equation}
 \vert \sum_{m=1}^{M}e^{-i\kappa\hat{x}\cdot z_m}Q_m\vert \leq M\sup_{m=1, ..., M}\vert Q_m\vert=O({\textcolor{black}{\epsilon}}^{1-s}). 
\end{equation}
Hence if the number of obstacles is $M:=M(\textcolor{black}{\epsilon}):=O({\textcolor{black}{\epsilon}}^{-s}), \; s<1$ and $t$ satisfies (\ref{general-condition-s-t}), ${\textcolor{black}{\epsilon}}\rightarrow 0$, 
then from (\ref{x oustdie1 D_m farmainp-near*}, \ref{x oustdie1 D_m farmains-near*}), we deduce that
\begin{equation}\label{s-smaller-1}
 U^\infty(\hat{x},\theta)\rightarrow 0, \mbox{ as } {\textcolor{black}{\epsilon}}\rightarrow 0, \mbox{ uniformly in terms of } \theta \mbox{ and } \hat{x} \mbox{ in } \mathbb{S}^2.
\end{equation}
This means that this collection of obstacles has no effect on the homogeneous medium as ${\textcolor{black}{\epsilon}}\rightarrow 0$. 
\bigskip

\begin{figure}
\centering
 \input{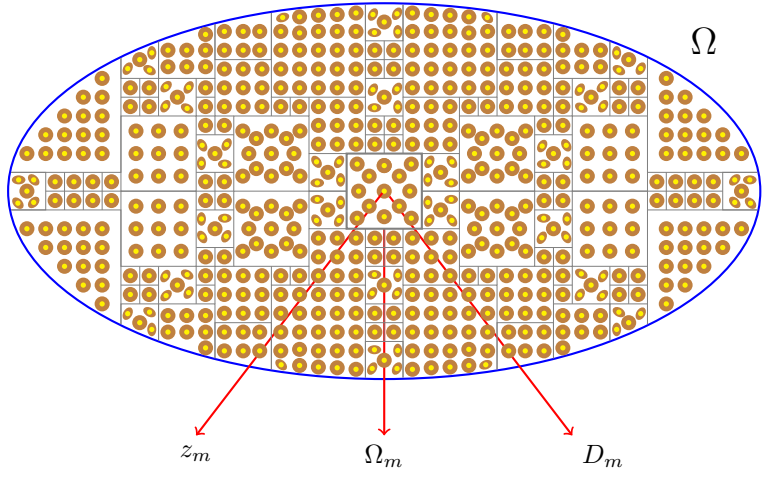}
\caption{A Schematic example on how the obstacles are distributed in $\Omega$.}\label{distribution-obstacles}
\end{figure}

Let us consider the case when $s=1$. We set $\Omega$ to be a bounded domain, say of unit volume, containing the obstacles $D_m, m=1, ..., M$. 
Given a  positive and continuous function $K: \mathbb{R}^3\rightarrow \mathbb{R}$, we divide $\Omega$ into $[{\textcolor{black}{\epsilon}}^{-1}]$ subdomains $\Omega_m,\; m=1, ..., [{\textcolor{black}{\epsilon}}^{-1}]$, each of volume $\textcolor{black}{\epsilon}\frac{[K(z_m)+1]}{K(z_m)+1}$, with $z_m \in \Omega_m$ as its center and contains $[K(z_m)+1]$ obstacles, see Fig \ref{distribution-obstacles}. 
We set $K_{max}:=\sup_{z_m}(K(z_m)+1)$, hence $M=\sum^{[{\textcolor{black}{\epsilon}}^{-1}]}_{j=1}[K(z_m)+1]\leq K_{max}[{\textcolor{black}{\epsilon}}^{-1}]=O({\textcolor{black}{\epsilon}}^{-1})$. 

\begin{theorem}\label{equivalent-medimu}
Let the small obstacles be distributed in a bounded domain $\Omega$, say of unit volume, with their number $M:=M(\textcolor{black}{\epsilon}):=O(\textcolor{black}{\epsilon}^{-1})$ and their minimum distance 
$d:=d(\textcolor{black}{\epsilon}):=O(\textcolor{black}{\epsilon}^{t})$,\; $\frac{1}{3}\leq t < \frac{1}{2}$, as $\textcolor{black}{\epsilon}\rightarrow 0$, as described above. In addition, we assume that the Lam\'e coefficients $\lambda$ and $\mu$ satisfy the conditions $
\max\{{\kappa_{s}},\; {\kappa_{p}}\} < \frac{2}{\diam(\Omega)},$ and $
\mathring{c}(\lambda +2 \mu) <\frac{\pi}{\sqrt{14}M_{max}\,{\max_m}C^a(B_m)}
$ \footnote{The constant $\mathring{c}$ is defined in (\ref{constant-final-estimates-Green}) and  ${C}^a_m$ denotes the capacitance of each scatterer (i.e. the acoustic capacitance). These conditions on the Lam\'e coefficients $\lambda$ and $\mu$ can be replaced by a condition on the wave number $\omega$, see Lemma \ref{Mazyawrkthmela} and the related footnote.}.
\begin{enumerate}
 
 \item If the obstacles are distributed arbitrarily in $\Omega$, i.e. with different capacitances, then there exists a potential $\bold{C}_0 \in \cap_{p\geq 1}L^p(\mathbb{R}^{3})$ with support in $\Omega$ such that
 \begin{equation}\label{B}
\lim_{a\rightarrow 0}U^\infty(\hat{x},\theta)= U_{0}^\infty(\hat{x},\theta) \mbox{ uniformly in terms of } \theta \mbox{ and } \hat{x} \mbox{ in } \mathbb{S}^2
 \end{equation}
where $U_{0}^\infty(\hat{x},\theta)$ is the farfield corresponding to the scattering problem

\begin{equation}
(\Delta^e + \omega^{2}-(K+1)\bold{C}_0)U_{0}^{t}=0 \mbox{ in }\mathbb{R}^{3},
\end{equation}

\begin{equation}
U_0^{t}|_{\partial D_m}=0,\, 1\leq m \leq M
\end{equation}

with the radiation conditions
\begin{equation}
\lim_{|x|\rightarrow\infty}|x|(\frac{\partial U_{0, p}}{\partial|x|}-i \textcolor{black}{\kappa_{p}}U_{0, p})=0,  \mbox{  and  }
\lim_{|x|\rightarrow\infty}|x|(\frac{\partial U_{0, s}}{\partial|x|}-i \textcolor{black}{\kappa_{s}}U_{0, s})=0
\end{equation}

\item If in addition $K\mid_{\Omega}$ is in $C^{0, \gamma}(\Omega)$, $\gamma \in (0, 1]$ and the obstacles have the same capacitances $C$, then
\begin{equation}\label{C}
  U^\infty(\hat{x},\theta)= U_0^\infty(\hat{x},\theta) +O(\textcolor{black}{\epsilon}^{\min\{\gamma, \frac{1}{3}, \frac{3}{2}-3t\}}) \mbox{ uniformly in terms of } \theta \mbox{ and } \hat{x} \mbox{ in } \mathbb{S}^2
 \end{equation}
 where $\bold{C}_0=C$ in $\Omega$ and $\bold{C}_0=0$ in $\mathbb{R}^{3} \setminus{\overline \Omega}$.
\end{enumerate}

\end{theorem}

\subsection{The results for variable background elastic mass density}

\bigskip

Assume that the Lam\'e coefficients $\lambda$ and $\mu$ are constants satisfying $ \mu > 0 \mbox{ and } 3\lambda+2\mu >0$ and the mass density $\rho$ to be a measurable and bounded function which is equal to a constant that we normalize to a unity outside of a bounded domain $\Omega$. We set $\rho_{max}$ to be the upper bound of $\rho$.

In this case, the total field $U_\rho^{t}:=U^{i}+U_\rho^{s}$
satisfies the following exterior Dirichlet problem of the elastic waves
\begin{equation}
(\Delta^e + \omega^{2}\rho)U_\rho^{t}=0 \mbox{ in }\mathbb{R}^{3}\backslash \left(\mathop{\cup}_{m=1}^M \bar{D}_m\right),\label{elaimpoenetrable-rho}
\end{equation}
\begin{equation}
U_\rho^{t}|_{\partial D_m}=0,\, 1\leq m \leq M\label{elagoverningsupport-rho}
\end{equation}

with the Kupradze radiation conditions (K.R.C)
\begin{equation}\label{radiationcela-rho}
\lim_{|x|\rightarrow\infty}|x|(\frac{\partial U_{\rho, p}}{\partial|x|}-i \textcolor{black}{\kappa_{p}}U_{\rho, p})=0,  \mbox{  and  }
\lim_{|x|\rightarrow\infty}|x|(\frac{\partial U_{\rho, s}}{\partial|x|}-i \textcolor{black}{\kappa_{s }}U_{\rho, s})=0,
\end{equation}
where the two limits are uniform in all the directions $\hat{x}:=\frac{x}{|x|}\in \mathbb{S}^{2}$ and $U_{\rho, p}$ and $U_{\rho, s}$ are respectively the P-part and S-part of the scattered field $U^s_{\rho, p}$

The scattering problem (\ref{elaimpoenetrable-rho}-\ref{radiationcela-rho}) is well posed in the H\"{o}lder or Sobolev spaces,
see \cite{C-K:1998,C-K:1983,Kupradze:1965, K-G-B-B:1979} for instance, and the scattered field $U^s$ has the following asymptotic expansion:

\begin{equation}\label{Lamesystemtotalfieldasymptoticsmall-rho}
 U_\rho^s(x) := \frac{e^{i\textcolor{black}{\kappa_{p}}|x|}}{|x|}U^{\infty}_{\rho, p}(\hat{x}) +
\frac{e^{i\textcolor{black}{\kappa_{s}}|x|}}{|x|}U^{\infty}_{\rho, s}(\hat{x}) + O(\frac{1}{|x|^{2}}),~|x|\rightarrow \infty
\end{equation}
uniformly in all directions $\hat{x}\in \mathbb{S}^{2}$. The longitudinal part of the far-field, i.e. $U^{\infty}_{\rho, p}(\hat{x})$ is normal to $\mathbb{S}^{2}$
while the transversal part $U^{\infty}_{s}(\hat{x})$ is tangential to $\mathbb{S}^{2}$. We set $U_\rho^{\infty}:=(U_{\rho, p}^{\infty}, U_{\rho, s}^{\infty})$.

As in the case of constant background mass density, there exist two positive constants $\textcolor{black}{\epsilon}_0$ and $c_0$ depending only on the size of $\Omega$, the
Lipschitz character of $B_m,m=1,\dots,M$, $d_{\max}$, $\omega_{\max}$, $\rho_{max}$ and $\rho_{max}$ such that
if
\begin{equation}\label{conditions-elasma-rho}
\textcolor{black}{\epsilon} \leq \textcolor{black}{\epsilon}_0 ~~ \mbox{and} ~~ \sqrt{M-1}\frac{\textcolor{black}{\epsilon}}{d}\leq c_0
\end{equation}
 then we have the following asymptotic expansion
 for the P-part, $U^\infty_{\rho, p}(\hat{x},\theta)$, and the S-part, $U^\infty_{\rho, s}(\hat{x},\theta)$, of the far-field pattern:
 \begin{eqnarray}
  &\hspace{-1cm}U^\infty_{\rho, p}(\hat{x},\theta)=& V^\infty_{\rho, p}(\hat{x},\theta)\;+~ \left[\sum_{m=1}^{M}G_{\rho, p}^{\infty}(\hat{x}, z_m)Q_{\rho, m}\right.\nonumber\\
  &&\hspace{-2cm}\left.+O\left(M\left[\textcolor{black}{\epsilon}^2+\frac{\textcolor{black}{\epsilon}^3}{d^{5-3\alpha}}+\frac{ \textcolor{black}{\epsilon}^4}{d^{9-6\alpha}}\right]+M(M-1)\left[\frac{\textcolor{black}{\epsilon}^3}{d^{2\alpha}}+\frac{\textcolor{black}{\epsilon}^4}{d^{4-\alpha}}+
  \frac{\textcolor{black}{\epsilon}^4}{d^{5-2\alpha}}\right]+M(M-1)^2\frac{\textcolor{black}{\epsilon}^4}{d^{3\alpha}}\right) \right],
   \label{x oustdie1 D_m farmainp-near-rho}\\
 &\hspace{-1cm}U^\infty_{\rho, s}(\hat{x},\theta)=& V^\infty_{\rho, s}(\hat{x},\theta)\;+~ \left[\sum_{m=1}^{M}G_{\rho, s}^{\infty}(\hat{x}, z_m) Q_{\rho, m}\right.\nonumber\\
 &&\hspace{-2cm}\left.+O\left(M\left[\textcolor{black}{\epsilon}^2+\frac{\textcolor{black}{\epsilon}^3}{d^{5-3\alpha}}+\frac{ \textcolor{black}{\epsilon}^4}{d^{9-6\alpha}}\right]+M(M-1)\left[\frac{\textcolor{black}{\epsilon}^3}{d^{2\alpha}}+\frac{\textcolor{black}{\epsilon}^4}{d^{4-\alpha}}+
 \frac{\textcolor{black}{\epsilon}^4}{d^{5-2\alpha}}\right]+M(M-1)^2\frac{\textcolor{black}{\epsilon}^4}{d^{3\alpha}}\right) \right],
 \label{x oustdie1 D_m farmains-near-rho}
  \end{eqnarray}
where $G_{\rho, p}^{\infty}(\hat{x}, z_m)$ and $G_{\rho, s}^{\infty}(\hat{x}, z_m)$ are the P-part and S-part of the farfields of the Green's function $G_\rho(x, z)$, of the operator $\Delta^e + \omega^{2}\rho$ in the whole space $\mathbb{R}^3$, evaluated in the direction $\hat{x}$ and the source point $z_m$.
 
 The vector coefficients $Q_{\rho, m}$, $m=1,..., M,$ are the solutions of the following linear algebraic system
\begin{eqnarray}\label{fracqcfracmain-rho}
 C_m^{-1}Q_{\rho, m} &=&-V_\rho(z_m, \theta)-\sum_{\substack{j=1 \\ j\neq m}}^{M} G_\rho(z_m,z_j)Q_{\rho, j},~~
\end{eqnarray}
for $ m=1,..., M,$ with $V_\rho(\cdot, \theta):=V_\rho^s(\cdot, \theta)+U^i(\cdot,\theta)$ is the total field satisfying 
\begin{equation}
(\Delta^e + \omega^{2}\rho)V_\rho^{t}=0 \mbox{ in }\mathbb{R}^{3},\label{elaimpoenetrable-rho1}
\end{equation}
and the scattered field $V_\rho^s(\cdot, \theta)$ the Kupradze radiation conditions (K.R.C).
 
\bigskip

\begin{corollary}\label{equivalent-medimu-rho}
Let the small obstacles be distributed in a bounded domain $\Omega$, say of unit volume, with their number $M:=M(\textcolor{black}{\epsilon}):=O(\textcolor{black}{\epsilon}^{-1})$ and their minimum distance 
$d:=d(\textcolor{black}{\epsilon}):=O(\textcolor{black}{\epsilon}^{t})$,\; $\frac{1}{3}\leq t < \frac{1}{2}$, as $\textcolor{black}{\epsilon}\rightarrow 0$, as described above. In addition, we assume that the Lam\'e coefficients $\lambda$ and $\mu$ satisfy the conditions $
\max\{{\kappa_{s}},\; {\kappa_{p}}\} < \frac{2}{\diam(\Omega)},$ and $
\mathring{c}(\lambda +2 \mu) <\frac{\pi}{\sqrt{14}M_{max}\,\textcolor{black}{{{\max_m}C^a(B_m)}}}
$. \footnote{Same comments on $\mathring{c}$ as in Theorem \ref{equivalent-medimu} apply here too.}
\begin{enumerate}
 
 \item If the obstacles are distributed arbitrarily in $\Omega$, i.e. with different capacitances, then there exists a potential $\bold{C}_0 \in \cap_{p\geq 1}L^p(\mathbb{R}^{3})$ with support in $\Omega$ such that
 \begin{equation}\label{B-rho}
\lim_{a\rightarrow 0}U_\rho^\infty(\hat{x},\theta)= U_{\rho, 0}^\infty(\hat{x},\theta) \mbox{ uniformly in terms of } \theta \mbox{ and } \hat{x} \mbox{ in } \mathbb{S}^2
 \end{equation}
where $U_{\rho, 0}^\infty(\hat{x},\theta)$ is the farfield corresponding to the scattering problem

\begin{equation}
(\Delta^e + \omega^{2}\rho-(K+1)\bold{C}_0)U_{_\rho, 0}^{t}=0 \mbox{ in }\mathbb{R}^{3},
\end{equation}

\begin{equation}
U_{\rho, 0}^{t}|_{\partial D_m}=0,\, 1\leq m \leq M
\end{equation}

with the radiation conditions.

\item If in addition $K\mid_{\Omega}$ is in $C^{0, \gamma}(\Omega)$, $\gamma \in (0, 1]$ and the obstacles have the same capacitances, then
\begin{equation}\label{C-rho}
  U_\rho^\infty(\hat{x},\theta)= U_{\rho, 0}^\infty(\hat{x},\theta) +O(\textcolor{black}{\epsilon}^{\min\{\gamma, \frac{1}{3}, \frac{3}{2}-3t\}}) \mbox{ uniformly in terms of } \theta \mbox{ and } \hat{x} \mbox{ in } \mathbb{S}^2
 \end{equation}
 where $\bold{C}_0=C$ in $\Omega$ and $\bold{C}_0=0$ in $\mathbb{R}^{3} \setminus{\overline \Omega}$.
\end{enumerate}

\end{corollary}

\subsection{Applications of the results and a comparison to the literature}
The main contribution of this work is to have shown that by removing from a bounded region of an elastic background, modeled by constant Lam\'e coefficients $\lambda$ and $\mu$ 
and a possibly variable density $\rho$, a number $M:=M(\textcolor{black}{\epsilon})\sim \textcolor{black}{\epsilon}^{-1}$  of small and rigid bodies with diameter of order $\textcolor{black}{\epsilon}$, 
distant from each other of at least 
$d:=d(\textcolor{black}{\epsilon})\sim \textcolor{black}{\epsilon}^t$, $\frac{1}{3}\leq t < \frac{1}{2}$, then the 'perforated' medium behaves, as $\textcolor{black}{\epsilon} \rightarrow 0$, as a new elastic medium modeled by the same 
Lam\'e coefficients $\lambda$ and $\mu$ but with a new coefficient $\omega^2\rho-(K +1)\bold{C}_0$. The coefficient $K$ models the local distribution (or the local number) 
of the bodies while the coefficient $\bold{C}_0$, coming from the capacitance of the bodies, describes the geometry of the small bodies as well as their elastic directional 
diffusion properties (i.e. the anisotropy character). In addition, we provide explicit error estimates between the far-fields corresponding to the perforated medium and 
the equivalent one. From this result we can make the following conclusions:
\begin{enumerate}
\item Assume that the removed bodies have spherical shapes. For these shapes the corresponding elastic capacitance $C$ is of the form $c I_3$ (i.e. a scalar multiplied by the identity matrix). In section \ref{structure-capacitance}, we describe a more general set of shapes satisfying this property. Hence the equivalent coefficient $\omega^2\rho-(K +1)c$ is isotropic while for general shapes it might be anisotropic. To achieve anisotropic coefficients, a possible choice of the shapes might be an ellipse.   

\item If we choose the local number of bodies $K$ large enough or the shapes of the reference bodies, $B_m$, $m=1,...,M$, having a large capacitance (i.e. a relative large radius) so that $\omega^2\rho-(K +1)c<0$. This means that, at the frequencies $\omega$ for which this inequality is satisfied, the elastic material, located in $\Omega$, should have an unusual behavior, namely the traction forces (modeled by the Lam\'e coefficients $\lambda$ and $\mu$) and the addition of the two body forces (modeled by the density $\rho$ and the coefficient $(K +1)\bold{C}_0$ respectively) will act in the same direction.

\item Assume that the background medium is modeled by variable mass density $\rho>1$ in $\Omega$. If we remove small bodies from $\Omega$ with appropriate $K$ and/ or capacitance $\bold{C}_0$ so that $\rho-(K +1)c\omega^{-2}=1$, then the new elastic material will behave every where in $\mathbb{R}^3$ as the background medium. Hence the new material will not scatter the sent incident waves at the given frequency $\omega$, i.e. the region $\Omega$ modeled by $\rho$ will be cloaked at that frequency.   
\end{enumerate}     

The 'equivalent' behavior between a collection of, appropriately dense, small holes
and an extended penetrable obstacle modeled by an additive potential was already observed by Cioranescu and Murat \cite{DC-FM:Book:BB1979, DC-FM:Book:BB1997} and also the references therein, where the coefficient $K$ is reduced to zero since locally they have only one hole. Their analysis is based on the homogenization theory for which they assume that the obstacles are distributed periodically, see also \cite{B-L-P:1978} and \cite{J-K-O:1994}. 

In the results presented here, we do not need such periodicity and no homogenization is used. Instead, the analysis is based on the invertibility properties of the algebraic system (\ref{fracqcfracmain}) 
and the precise treatment of the summation in the dominant terms of (\ref{x oustdie1 D_m farmainp-near})-(\ref{x oustdie1 D_m farmains-near}). This analysis was already tested for 
the acoustic model in \cite{DPC-SM13-3}. Compared to \cite{DPC-SM13-3}, here, in addition to the difficulties coming from the vector character of the Lam\'e system, we improved the order of the error estimate, 
i.e. $O(\textcolor{black}{\epsilon}^{\min\{\gamma, \frac{1}{3}, \frac{3}{2}-3t\}})$ instead of $O(\textcolor{black}{\epsilon}^{\min\{\gamma, \frac{1}{3}-\frac{4}{5}t\}})$ which, 
for $t:=\frac{1}{3}$ for instance, reduce to $O(\textcolor{black}{\epsilon}^{\min\{\gamma, \frac{1}{3}\}})$ and $O(\textcolor{black}{\epsilon}^{\min\{\gamma, \frac{1}{15}}\})$ respectively. 

Let us mention that a result similar to (\ref{B-rho}), for the acoustic model, is also derived by Ramm in several of his papers, see for instance \cite{RAMM:2007}, but without error estimates. Compared to his results, and as we said earlier in addition to the vector character of Lame model, we provide the approximation by improved explicit error estimates without any other assumptions while,  as shown in \cite{RAMM:2007, RAMM:2011} for instance, in addition to some formal arguments, he needs extra 
assumptions on the distribution of the obstacles. 

We discussed above the application of our result in elastic cloaking via perforation (or perturbation with small bodies). Another way of achieving elastic cloaking is by using transformations, see \cite{HG-LH:JMPANS2015}.

The rest of the paper is organized as follows. In section 2, we give the detailed proof of Theorem \ref{equivalent-medimu}.  In section 3, we describe the one of Corollary \ref{equivalent-medimu-rho}
 by discussing the main changes one needs to make in the proof of Theorem \ref{equivalent-medimu}. Finally, in section 4, we discuss some invariant properties of the elastic capacitance to characterize the shapes that have a 'scalar' capacitance.

\section{Proof of Theorem \ref{equivalent-medimu}}

\subsection{The fundamental solution}\label{fdelsms}
The Kupradze matrix $\Gamma^\omega=(\Gamma^\omega_{ij})^3_{i,j=1}$ of the fundamental solution to the Navier equation is given by
\begin{eqnarray}\label{kupradzeten}
 \Gamma^\omega(x,y)=\frac{1}{\mu}\Phi_{\textcolor{black}{\kappa_{s}}}(x,y)\rm \textbf{I}+\frac{1}{\omega^2}\nabla_x\nabla_x^{\top}[\Phi_{\textcolor{black}{\kappa_{s}}}(x,y)-\Phi_{\textcolor{black}{\kappa_{p}}}(x,y)],
\end{eqnarray}
where $\Phi_{\kappa}(x,y)=\frac{e^{i\kappa|x-y|}}{4\pi|x-y|}$ denotes the free space fundamental solution of the Helmholtz equation $(\Delta+\kappa^2)\,u=0$ in
$\mathbb{R}^3$. The asymptotic behavior of Kupradze tensor at infinity is given as follows
\begin{equation}\label{elafundatensorasymptotic}
 \Gamma^{\omega}(x,y)=\frac{1}{4\pi\,c_p^{2}}\hat{x}\otimes\hat{x} \frac{e^{i\textcolor{black}{\kappa_{p}}|x|}}{|x|}e^{-i\textcolor{black}{\kappa_{p}}\hat{x}\cdot\,y} +
\frac{1}{4\pi\,c_s^{2}}({\rm \textbf{I}}- \hat{x}\otimes\hat{x}) \frac{e^{i\textcolor{black}{\kappa_{s}}|x|}}{|x|}e^{-i\textcolor{black}{\kappa_{s}}\hat{x}\cdot\,y}+O(|x|^{-2})
\end{equation}
with $\hat{x}=\frac{x}{|x|}\in\mathbb{S}^{2}$ and $\rm \textbf{I}$ being the identity matrix in $\mathbb{R}^{3}$, see \cite{A-K:IMA2002} for instance.
As mentioned in \cite{A-K-L-LameCPDE2007},  \eqref{kupradzeten} can also be represented as
\begin{eqnarray}\label{kupradzeten1}
 \Gamma^\omega(x,y)&=&\frac{1}{4\pi}\sum_{l=0}^{\infty}\frac{i^l}{l!(l+2)}\frac{1}{\omega^2}\left((l+1)\textcolor{black}{\kappa_{s}^{l+2}}+\textcolor{black}{\kappa_{p}^{l+2}}\right)|x-y|^{l-1}\rm \textbf{I}\nonumber\\
         & &-\frac{1}{4\pi}\sum_{l=0}^{\infty}\frac{i^l}{l!(l+2)}\frac{(l-1)}{\omega^2}\left(\textcolor{black}{\kappa_{s}^{l+2}}-\textcolor{black}{\kappa_{p}^{l+2}}\right)|x-y|^{l-3}(x-y)\otimes(x-y),
               \end{eqnarray}
from which we can get the gradient 
\begin{eqnarray}\label{gradkupradzeten1}
\hspace{-.5cm}\nabla_y \Gamma^\omega(x,y)&=&-\frac{1}{4\pi}\sum_{l=0}^{\infty}\frac{i^l}{l!(l+2)}\frac{(l-1)}{\omega^2}\left[\left((l+1)\textcolor{black}{\kappa_{s}^{l+2}}+\textcolor{black}{\kappa_{p}^{l+2}}\right)|x-y|^{l-3}(x-y)\textcolor{black}{\otimes}\rm \textbf{I}\right.\nonumber\\
&&-\left.\left(\textcolor{black}{\kappa_{s}^{l+2}}-\textcolor{black}{\kappa_{p}^{l+2}}\right)|x-y|^{l-3}\left((l-3)|x-y|^{-2}\otimes^{3}(x-y)+\rm \textbf{I}\textcolor{black}{\otimes}(x-y)+\sum_{j=1}^3\sum_{j^\prime=1}^3(x-y)_j\,e_{j^\prime}\textcolor{black}{\otimes}\,e_{j}\textcolor{black}{\otimes}\,e_{j^\prime}\right)\right],\nonumber\\
\end{eqnarray}
where $e_j$ is standard basis vector. 
Using the formulas \eqref{kupradzeten1} and \eqref{gradkupradzeten1} 
we can have the following estimates, for $x,y\in\Omega, x\neq\;y$, see \cite{C-S:2015};
 \begin{eqnarray}
\left|\Gamma^\omega(x,y)\right|
\leq\frac{1}{4\pi}\left[\frac{C_7}{\vert\;x-y\vert}+C_8\right],\qquad&&
\left|\nabla_y\Gamma^\omega(x,y)\right|
\leq  \frac{1}{4\pi}\left[\frac{C_9}{\vert\;x-y\vert^2}+C_{10}\right],\label{gradkupradzeten1dfnt}
\end{eqnarray}
with

 $$C_7:=\left[\frac{{1}}{c_s^2}+\frac{\textcolor{black}{1}}{c_p^2}\right],\qquad C_9:=\textcolor{black}{7}\left(\frac{1}{c_s^2}+\frac{1}{c_p^2}\right),$$
$$C_8:=\textcolor{black}{2}\frac{\omega}{c_s^3}\left(\frac{1-\left(\frac{1}{2}\textcolor{black}{\kappa_{s}}\diam(\Omega)\right)^{N_\Omega}}{1-\frac{1}{2}\textcolor{black}{\kappa_{s}}\diam(\Omega)}
+\frac{1}{2^{N_{\Omega}-1}}\right)+\frac{\omega}{c_p^3}\left(\frac{1-\left(\frac{1}{2}\textcolor{black}{\kappa_{p}}\diam(\Omega)\right)^{N_\Omega}}{1-\frac{1}{2}\textcolor{black}{\kappa_{p}}\diam(\Omega)}+\frac{1}{2^{N_{\Omega}-1}}\right),
$$ 
$$C_{10}:=\textcolor{black}{6}\frac{\omega^2}{c_s^4}\left(\textcolor{black}{\frac{1}{24}}+\frac{1-\left(\frac{1}{2}\textcolor{black}{\kappa_{s}}\diam(\Omega)\right)^{N_\Omega}}{1-\frac{1}{2}\textcolor{black}{\kappa_{s}}\diam(\Omega)}+\frac{1}{2^{N_{\Omega}-1}}\right)
+\textcolor{black}{5}\frac{\omega^2}{c_p^4}\left(\textcolor{black}{\frac{1}{20}}+\frac{1-\left(\frac{1}{2}\textcolor{black}{\kappa_{p}}\diam(\Omega)\right)^{N_\Omega}}{1-\frac{1}{2}\textcolor{black}{\kappa_{p}}\diam(\Omega)}+\frac{1}{2^{N_{\Omega}-1}}\right),
$$
and $N_{\Omega}=[2\diam(\Omega)\max\{\textcolor{black}{\kappa_{s}},\textcolor{black}{\kappa_{p}}\}e^2]$ where we assume the frequency  $\omega$ and the Lam\'e parameters $\lambda$ and $\mu$  to satisfy the condition $\max\{\textcolor{black}{\kappa_{s}},\; \textcolor{black}{\kappa_{p}}\} < \frac{2}{\diam(\Omega)}$.

\bigskip

\begin{indeed}
 for $x\in \bar{D}_m$ and $s\in \bar{D}_j$, first let us recall the following result from \cite[Lemma 2.6]{C-S:2015}.
 \begin{lemma}\label{stirlingapproxlemma}
 For each $k>0$ and for every $n\in\mathbb{Z}^{+}$ with $n\geq{ke^2}\, [=:N(k)]$ we have $n!\geq k^{n-1}$.
\end{lemma}
Now, by making use of \eqref{kupradzeten1}, we get the following estimate;
 \begin{eqnarray*}\label{modgamma_{ij}ela}
\left|\Gamma^\omega(x,y)\right|
&\leq& \frac{1}{4\pi}\frac{1}{\omega^2}\left(\textcolor{black}{\kappa_{s}^{2}}+\textcolor{black}{\kappa_{p}^{2}}\right)|x-y|^{-1}
\nonumber\\   &&
   +\frac{1}{4\pi}\sum_{l=1}^{\infty}\frac{1}{(l-1)!(l+2)}\frac{1}{\omega^2}
    \left(2\textcolor{black}{\kappa_{s}^{l+2}}+\textcolor{black}{\kappa_{p}^{l+2}}\right)|x-y|^{l-1} \nonumber\\
&\leq&  \frac{1}{4\pi}\frac{1}{\omega^2}\left[\frac{1}{|x-y|}\left(\textcolor{black}{\kappa_{s}^{2}}+\textcolor{black}{\kappa_{p}^{2}}\right)\right.
   \left.+\sum_{l=1}^{\infty}\frac{1}{(l-1)!(l+2)}\left(2\textcolor{black}{\kappa_{s}^{l+2}}+\textcolor{black}{\kappa_{p}^{l+2}}\right)\diam(\Omega)^{l-1}\right]\nonumber\\
&\leq&  \frac{1}{4\pi}\left[\frac{1}{|x-y|}\left(\frac{1}{c_s^2}+\frac{1}{c_p^2}\right)\right.
  \left.+\sum_{l=1}^{\infty}\frac{1}{(l-1)!(l+2)}\left(2\frac{\omega}{c_s^3}\textcolor{black}{\kappa_{s}^{l-1}}+\frac{\omega}{c_p^3}\textcolor{black}{\kappa_{p}^{l-1}}\right)\diam(\Omega)^{l-1}\right] \nonumber\\
&&\mbox{[By writing $N_{\Omega}=[2\diam(\Omega)\max\{\textcolor{black}{\kappa_{s}},\textcolor{black}{\kappa_{p}}\}e^2]$ and using Lemma \ref{stirlingapproxlemma}]}\nonumber\\
&\leq&  \frac{1}{4\pi}\left[\frac{1}{|x-y|}\left(\frac{1}{c_s^2}+\frac{1}{c_p^2}\right)\right.
 \left.+2\frac{\omega}{c_s^3}\left(\sum_{l=1}^{N_{\Omega}}\left(\frac{1}{2}\textcolor{black}{\kappa_{s}}\diam(\Omega)\right)^{l-1}+\sum_{l=N_{\Omega}+1}^{\infty}\frac{1}{2^{l-1}}\right)\right.\nonumber\\
&&\hspace{5cm}\left.+\frac{\omega}{c_p^3}\left(\sum_{l=1}^{N_{\Omega}}\left(\frac{1}{2}\textcolor{black}{\kappa_{p}}\diam(\Omega)\right)^{l-1}+\sum_{l=N_{\Omega}+1}^{\infty}\frac{1}{2^{l-1}}\right)\right] \nonumber\\
&=&  \frac{1}{4\pi}\left[\frac{1}{|x-y|}\left(\frac{1}{c_s^2}+\frac{1}{c_p^2}\right)\right.
 \left.+2\frac{\textcolor{black}{\kappa_{s}}}{c_s^2}\left(\frac{1-\left(\frac{1}{2}\textcolor{black}{\kappa_{s}}\diam(\Omega)\right)^{N_\Omega}}{1-\frac{1}{2}\textcolor{black}{\kappa_{s}}\diam(\Omega)}+\frac{1}{2^{N_{\Omega}-1}}\right)\right.\nonumber\\
&&\hspace{5cm}\left.+\frac{\textcolor{black}{\kappa_{p}}}{c_p^2}\left(\frac{1-\left(\frac{1}{2}\textcolor{black}{\kappa_{p}}\diam(\Omega)\right)^{N_\Omega}}{1-\frac{1}{2}\textcolor{black}{\kappa_{p}}\diam(\Omega)}+\frac{1}{2^{N_{\Omega}-1}}\right)\right]\nonumber\\
&=&\frac{1}{4\pi}\left[\frac{C_7}{|x-y|}+C_8\right]
\end{eqnarray*}

 Now, by making use of \eqref{gradkupradzeten1} we get the following estimate;

 \begin{eqnarray*}\label{gradkupradzeten1dfnt--}
\left|\nabla_x\Gamma^\omega(x,y)\right|
&\leq&\frac{1}{4\pi}\frac{1}{\omega^2}
 \left[7\left(\textcolor{black}{\kappa_{s}^{2}}+\textcolor{black}{\kappa_{p}^{2}}\right)|x-y|^{-2}+\frac{1}{8}\left(14\textcolor{black}{\kappa_{s}^{4}}+12\textcolor{black}{\kappa_{p}^{4}}\right)\right]  \nonumber\\
&&+\frac{1}{4\pi}\sum_{l=3}^{\infty}\frac{1}{l(l-2)!(l+2)}\frac{1}{\omega^2}
 \left((2l+8)\textcolor{black}{\kappa_{s}^{l+2}}+(l+8)\textcolor{black}{\kappa_{p}^{l+2}}\right)|x-y|^{l-2}  \nonumber\\
&\leq&\frac{1}{4\pi}\frac{1}{\omega^2}
 \left[\frac{7}{|x-y|^2}\left(\textcolor{black}{\kappa_{s}^{2}}+\textcolor{black}{\kappa_{p}^{2}}\right)+\frac{1}{4}\left(\textcolor{black}{\kappa_{s}^{4}}+\textcolor{black}{\kappa_{p}^{4}}\right)\right.\nonumber\\
&&\hspace{2cm}\left.+\sum_{l=2}^{\infty}\frac{1}{l(l-2)!(l+2)}
 \left((2l+8)\textcolor{black}{\kappa_{s}^{l+2}}+(l+8)\textcolor{black}{\kappa_{p}^{l+2}}\right)\diam(\Omega)^{l-2}\right]\nonumber\\
&\leq&  \frac{1}{4\pi}\left[\frac{7}{|x-y|^2}\left(\frac{1}{c_s^2}+\frac{1}{c_p^2}\right)
 +\frac{1}{4}\left(\frac{\omega^{2}}{c_s^4}+\frac{\omega^{2}}{c_p^4}\right)\right.\nonumber\\
&&\hspace{2cm}\left.+\sum_{l=2}^{\infty}\frac{1}{(l-2)!(l+2)}\left(6\frac{\omega^2}{c_s^4}\textcolor{black}{\kappa_{s}^{l-2}}+5\frac{\omega^2}{c_p^4}\textcolor{black}{\kappa_{p}^{l-2}}\right)\diam(\Omega)^{l-2}\right] \nonumber\\
&&\mbox{[By recalling $N_{\Omega}=[2\diam(\Omega)\max\{\textcolor{black}{\kappa_{s}},\textcolor{black}{\kappa_{p}}\}e^2]$ and using Lemma \ref{stirlingapproxlemma}]}\nonumber\\
&\leq&  \frac{1}{4\pi}\left[\frac{7}{|x-y|^2}\left(\frac{1}{c_s^2}+\frac{1}{c_p^2}\right)
  +\frac{1}{4}\left(\frac{\omega^{2}}{c_s^4}+\frac{\omega^{2}}{c_p^4}\right)\right.\nonumber\\
&&\hspace{2cm}\left.+6\frac{\omega^2}{c_s^4}\left(\sum_{l=0}^{N_{\Omega}-1}\left(\frac{1}{2}\textcolor{black}{\kappa_{s}}\diam(\Omega)\right)^{l}+\sum_{l=N_{\Omega}}^{\infty}\left(\frac{1}{2}\right)^{l}\right)\right.\nonumber\\
&&\hspace{2cm}\left.+5\frac{\omega^2}{c_p^4}\left(\sum_{l=0}^{N_{\Omega}-1}\left(\frac{1}{2}\textcolor{black}{\kappa_{p}}\diam(\Omega)\right)^{l}+\sum_{l=N_{\Omega}}^{\infty}\left(\frac{1}{2}\right)^{l}\right)\right] \nonumber\\
&=&  \frac{1}{4\pi}\left[\frac{C_9}{|x-y|^2}+C_{10}\right].
\end{eqnarray*}
\end{indeed}

The estimates in (\ref{gradkupradzeten1dfnt}) can be written as 
 \begin{eqnarray}
\left|\Gamma^\omega(x,y)\right|
\leq\frac{\mathring{c}}{4\pi\vert\;x-y\vert},\qquad&&
\left|\nabla_y\Gamma^\omega(x,y)\right|
\leq  \frac{\mathring{c}}{4\pi\vert\;x-y\vert^2},\label{gradkupradzeten1dfnt-1}
\end{eqnarray}
for different points $x,y\in \Omega$, where 
\begin{equation}\label{constant-final-estimates-Green}
\mathring{c}:=\max\{C_7, C_8\; \diam(\Omega), C_8, C_{10}\; \diam(\Omega)\}.
\end{equation}

\subsection{The relative distribution of the small bodies}\label{The-relative-distribution-of-the-small-bodies}
\begin{figure}[htp]
\centering
\includegraphics[width=4.5cm,height =4.5cm]{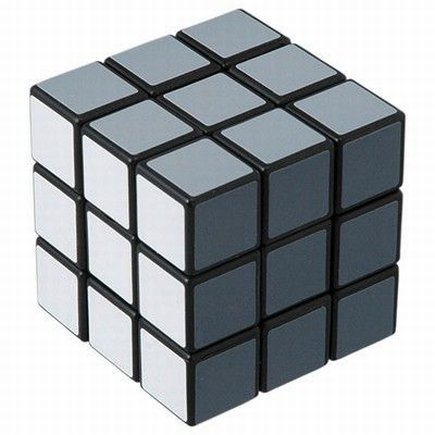}
\caption{Rubik's cube consisting of a {single} layer}\label{fig:1-acsmall}
\end{figure}

The following observation will be useful for the proof of Theorem \ref{equivalent-medimu}. 
For $m=1,\dots,M$ fixed, we distinguish between the obstacles $D_j$, $j\neq\,m$ by keeping them into different layers based on their distance from $D_m$.  
Let us first assume that $K(z_m)=0$ for every $z_m$. Hence each $\Omega_m$ has the (same) volume $\textcolor{black}{\epsilon}$ and contains only one obstacle $D_m$. 
Without loss of generality, we can take the $\Omega_m$'s as cubes. Hence we can suppose that these cubes are arranged in cuboids, 
for example Rubik's cube, in different layers and $\Omega_m$ is being located at the center. Observe that the total number of cubes in a Rubik's cube consisting of 
$N$ layers is $(2N+1)^3$ cubes, {see Fig \ref{fig:1-acsmall}}.  Since, the volume of each $\Omega_m$ is $\textcolor{black}{\epsilon}$  and the cubes are arranged in $N$ layered
Rubik's cube, the total number of cubes are $\textcolor{black}{\epsilon}^{-1}$ which is then equal to $(2N+1)^3$ and then $N=\left[\frac{\textcolor{black}{\epsilon}^{-\frac{1}{3}}-1}{2}\right]$.  
Hence, (1) the total cubes upto the $n^{th}$ layer consists of $(2n+1)^3$ cubes
and (2) the number of obstacles located in the $n^{th}$, $n\neq0$,  layer  is $[(2n+1)^3-(2n-1)^3]=24 n^2+2$, and their distance from $D_m$ is 
more than ${n}\left(\textcolor{black}{\epsilon}^{\frac{1}{3}}-\textcolor{black}{\epsilon}\frac{\textcolor{black}{\max\limits_{1\leq m\leq M } \diam (B_m)}}{2}\right)$, 
for $n=0,\dots,\left[\frac{\textcolor{black}{\epsilon}^{-\frac{1}{3}}-1}{2}\right]$.\footnote{Indeed, if we take any two adjacent cubes $\Omega_1$ and $\Omega_2$ which 
have the obstacles $D_1$ and $D_2$ at their centers respectively, then the distance between the obstacles  is more than 
$\left(\textcolor{black}{\epsilon}^{\frac{1}{3}}-\textcolor{black}{\epsilon}\frac{\textcolor{black}{\max\limits_{1\leq m\leq M } \diam (B_m)}}{2}\right)$ 
as the maximal radius of  $D_m$ is $\textcolor{black}{\epsilon}\frac{\textcolor{black}{\max\limits_{1\leq m\leq M } \diam (B_m)}}{2}$ and the volume of the cube $\Omega_m$ is $\textcolor{black}{\epsilon}$.} 
\par
Now, we come back to the case where $K(z_m)\neq 0$. First observe that $\frac{1}{2}\leq \frac{[K(z_m)+1]}{K(z_m)+1}\leq 1$. Hence with such $\Omega_m$'s, the total cubes located in 
the $n^{th}$ layer consists of at most the double of $[(2n+1)^3-(2n-1)^3]$, i.e. $48n^2+4$. It is due to the fact that upperbound of the fraction  $\frac{[K(z_m)+1]}{K(z_m)+1}$ 
\textcolor{black}{relates to} the lower bound of the total number of cubes ( which is explained in the previous paragraph for the case $\Omega_m$'s of volume $\textcolor{black}{\epsilon}$) 
and the lower bound of the fraction $ \frac{[K(z_m)+1]}{K(z_m)+1}$ \textcolor{black}{relates to} the upper bound of the total number of cubes located in each layer, which is $2[(2n+1)^3-(2n-1)^3]$.

\subsection{Solvability of the linear-algebraic system \eqref{fracqcfracmain}}\label{Solvability-of-the-linear-algebraic-system-elastic-small}

\bigskip

We start with the following lemma on the uniform bounds of the elastic capacitances, see \cite[Lemma 3.1]{C-S:2015} and the references therein.

\begin{lemma}\label{capacitance-eig-single} Let $\lambda^{min}_{eig_m}$ and $\lambda^{max}_{eig_m}$ be the minimal and maximal eigenvalues of the elastic capacitance matrices $\bar{C}_m$, for $m=1,2,\dots,M$. Denote
by ${C}^a_m$ the capacitance of each scatterer in the acoustic case,\footnote{Recall that, for $m=1,\dots,M$, ${C}^a_m:=\int_{\partial D_m}\sigma_m(s)ds$ and $\sigma_{m}$ is
the solution of the integral equation of the first kind $\int_{\partial D_m}\frac{\sigma_{m} (s)}{4\pi|t-s|}ds=1,~ t\in \partial D_m$, see \cite{DPC-SM13}.} then we have the following estimate;
\begin{eqnarray}\label{lowerupperestforintgradub-m}
 {\mu C^a(B_m)}\textcolor{black}{\epsilon} =\mu\,C^a_m\,\leq\,\lambda^{min}_{eig_m}\,\leq\,\lambda^{max}_{eig_m}\,\leq\,(\lambda+2\mu)\,C^a_m ={(\lambda +2 \mu)
 C^a(B_m)}\textcolor{black}{\epsilon},
\end{eqnarray}
for $ m=1,2,\dots,M$.
\end{lemma}
\bigskip

The constant $C(B_m)$ is the acoustic capacitance of the reference body $B_m$ which can be estimated above and below by the Lipschitz character of $B_m$, \cite{DPC-SM13}. The following lemma
provides us with the needed estimate on the invertibility of the algebraic system \eqref{fracqcfracmain} whose coefficient matrix '$\mathbf{B}$' is given by;
\begin{eqnarray}\label{coeffmatB}
\mathbf{B}:=\left(\begin{array}{ccccc}
   -\bar{C}_1^{-1} &-\Gamma^{\omega}(z_1,z_2)&-\Gamma^{\omega}(z_1,z_3)&\cdots&-\Gamma^{\omega}(z_1,z_M)\\
-\Gamma^{\omega}(z_2,z_1)&-\bar{C}_2^{-1}&-\Gamma^{\omega}(z_2,z_3)&\cdots&-\Gamma^{\omega}(z_2,z_M)\\
 \cdots&\cdots&\cdots&\cdots&\cdots\\
-\Gamma^{\omega}(z_M,z_1)&-\Gamma^{\omega}(z_M,z_2)&\cdots&-\Gamma^{\omega}(z_M,z_{M-1}) &-\bar{C}_M^{-1}
   \end{array}\right)
\end{eqnarray}

\begin{lemma}\label{Mazyawrkthmela}
The matrix $\mathbf{B}$ is invertible and the solution vector $Q_m\; m=1, ..., M,$ of \eqref{fracqcfracmain} satisfies the estimate:
\begin{eqnarray}\label{mazya-fnlinvert-small-ela-2}
 \sum_{m=1}^{M}\|Q_m\|_2^{2}
\leq \left( (\lambda +2 \mu)^{-1} - \frac{\sqrt{14}M_{max}\,\mathring{c}\,({\max_m} C^a_m)\;\textcolor{black}{\epsilon}^{-1}}{\pi} \right)^{-2}\left( C^a_m\right)^2\sum_{m=1}^{M}\|U^i(z_m)\|_2^2,
\end{eqnarray}
if $ {\max}C^a_m\; \textcolor{black}{\epsilon}^{-1}<\frac{\pi}{\sqrt{14}M_{max}\,\mathring{c}(\lambda +2 \mu)}$. In addition,
\begin{equation}\label{L1-norm-estimate-algebraic-system}
\sum_{m=1}^{M}\|Q_m\|_1
\leq \left( (\lambda +2 \mu)^{-1} - \frac{\sqrt{14}M_{max}\,\mathring{c}\,({\max_m} C^a_m)\;\textcolor{black}{\epsilon}^{-1}}{\pi} \right)^{-1}\;{\max}C^a_m\; M\; \max^M_{m =1}\|U^i(z_m)\|_2.
\end{equation}\footnote{
Alternatively, see \cite[Lemma 3.2]{C-S:2015}, we can have the following estimate (and then the invertibility of $\mathbf{B}$):

\begin{eqnarray*}\label{mazya-fnlinvert-small-ela-3}
\hspace{-.3cm} \sum_{m=1}^{M}\hspace{-.1cm}\|{Q}_m\|_2^{2}
\leq \hspace{-.07cm}4 \left(\min\limits_{1\leq m \leq M}\lambda^{min}_{eig_m}\hspace{-.07cm}-\hspace{-.07cm}\frac{3\tau}{5\pi\,d}\max\limits_{1\leq m \leq M}\lambda^{max^2}_{eig_m}\hspace{-.1cm}\right)^{\hspace{-.1cm}-2}\hspace{-.2cm}\left(\max\limits_{1\leq m \leq M}\lambda^{max}_{eig_m}\hspace{-.1cm}\right)^{\hspace{-.1cm}4}\sum_{m=1}^{M}\hspace{-.1cm}\|U^i(z_m)\|_2^2,
\end{eqnarray*}
if $\,\hspace{-.1cm}\left(\hspace{-.05cm}\max\limits_{1\leq\,m\leq\,M}\lambda^{max^2}_{eig_m}\hspace{-.05cm}\right)\hspace{-.1cm}<\hspace{-.1cm}{\tau}^{-1}\hspace{-.1cm}\left(\hspace{-.05cm}\frac{5\pi}{3}d\min\limits_{1\leq\,m\leq\,M}\lambda^{min}_{eig_m}\hspace{-.05cm}\right)\hspace{-.1cm}$ with the positively assumed value \\
\begin{center}
$
\tau\hspace{-.09cm}:=\hspace{-.1cm}\left[\hspace{-.05cm}\frac{1}{c_p^2}\hspace{-.1cm}-\hspace{-.1cm}2diam(\Omega)\frac{\omega}{c_s^3}\hspace{-.1cm}\left(\hspace{-.09cm}\frac{1-\left(\hspace{-.05cm}\frac{1}{2}\kappa_{s}diam(\Omega)\hspace{-.05cm}\right)^{N_\Omega}}{1-\left(\hspace{-.05cm}\frac{1}{2}\kappa_{s}diam(\Omega)\hspace{-.05cm}\right)}\hspace{-.1cm}+\hspace{-.1cm}\frac{1}{2^{N_{\Omega}-1}}\hspace{-.1cm}\right)
\hspace{-.1cm}-\hspace{-.1cm}diam(\Omega)\frac{\omega}{c_p^3}\hspace{-.1cm}\left(\hspace{-.09cm}\frac{1-\left(\hspace{-.05cm}\frac{1}{2}\kappa_{p}diam(\Omega)\hspace{-.05cm}\right)^{N_\Omega}}{1-\left(\hspace{-.05cm}\frac{1}{2}\kappa_{p}diam(\Omega)\hspace{-.05cm}\right)}\hspace{-.1cm}+\hspace{-.1cm}\frac{1}{2^{N_{\Omega}-1}}\hspace{-.1cm}\right)\hspace{-.1cm}\right]\hspace{-.1cm}
$
\end{center}
 with $\max\{{\kappa_{s}},\; {\kappa_{p}}\} < \frac{2}{\diam(\Omega)}$. Observe that this last condition and the positivity of $\tau$ can be seen as conditions on the wave number $\omega$. Using these conditions, we can replace in Theorem \ref{equivalent-medimu} and Corollary \ref{equivalent-medimu-rho} the conditions on the Lam\'e coefficients $\lambda$ and $\mu$ by a condition on $\omega$.
}
\end{lemma}
\bigskip
Since $C^a_m =C^a(B_m) \epsilon$
, then 
\begin{equation}\label{algebraic-system-invertible-condition}
{\max_m}C^a_m\; \textcolor{black}{\epsilon}^{-1}<\frac{\pi}{\sqrt{14}M_{max}\,\mathring{c}(\lambda +2 \mu)}
\end{equation} 
makes sense if $\label{algebraic-system-invertible-condition-2}
\textcolor{black}{{\max_m}C^a(B_m)}M_{max}<\frac{\pi}{\sqrt{14}\,\mathring{c}(\lambda +2 \mu)}.
$
As $C^a(B_m)$ is proportional to the radius of $B_m$, then (\ref{algebraic-system-invertible-condition}) 
will be satisfied if $\omega$ and the Lam\'e parameters $\lambda$ and $\mu$ satisfy the condition 
\begin{equation}\label{upper-bound-omega-lame-coe}
\max\{{\kappa_{s}},\; {\kappa_{p}}\} < \frac{2}{\diam(\Omega)} \mbox{ and }
\mathring{c}(\lambda +2 \mu) <\frac{\pi}{\sqrt{14}M_{max}\,}\textcolor{black}{\frac{1}{{\max_m}C^a(B_m)}}
\end{equation} 
recalling that $\mathring{c}$ is defined in (\ref{constant-final-estimates-Green}). Finally, let us observe that the right hand side of (\ref{upper-bound-omega-lame-coe}) depends only $M_{max}$ and the Lipschitz character of the reference obstacles $B_m$'s.

\begin{proofe}{\it{Proof of Lemma \ref{Mazyawrkthmela}.}} We start by factorizing $\mathbf{B}$ as $\mathbf{B}=-(\mathbf{C}^{-1}+\mathbf{B}_{n})$ where $\mathbf{C}:=Diag(\bar{C}_1,\bar{C}_2,\dots,\bar{C}_M)\in\mathbb{R}^{M\times\,M}$, $I$ is the identity matrix and 
$\mathbf{B}_{n}:=-\mathbf{C}^{-1}-\mathbf{B}$. We have $\mathbf{B}:\mathbb{C}^{3 M}\rightarrow\mathbb{C}^{3M}$, so it is enough to prove the injectivity in order to prove 
its invertibility. For this purpose,  let $X,Y$ are vectors in $\mathbb{C}^{3 M}$ and
 consider the system
\begin{eqnarray}\label{systemsolve1-small-ac}(\mathbf{C}^{-1}+\mathbf{B}_{n})X&=&Y.\end{eqnarray}
We denote by ${(\cdot)}^{real}$ and ${(\cdot)}^{img}$ the real and the imaginary parts of the corresponding complex vector/matrix. 
we also set $\mathbf{C}^{-1}$ by $\mathbf{C}_{I}$. From \eqref{systemsolve1-small-ac} we derive the following two identities:
\begin{eqnarray}
 (\mathbf{C}_{I}+\mathbf{B}^{real}_{n})X^{real}-\mathbf{B}^{img}_{n}X^{img}&=&Y^{real}\label{systemsolve1-small-sub1-ac-1},\\
 (\mathbf{C}_{I}+\mathbf{B}^{real}_{n})X^{img}+\mathbf{B}^{img}_{n}X^{real}&=&Y^{img}\label{systemsolve1-small-sub2-ac-1},
\end{eqnarray}
and then
\begin{eqnarray}
 \langle\,(\mathbf{C}_{I}+\mathbf{B}^{real}_{n})X^{real},X^{real}\rangle\,-\langle\,\mathbf{B}^{img}_{n}X^{img},X^{real}\rangle&=&\langle\,Y^{real},X^{real}\rangle\label{systemsolve1-small-sub1-ac-2},\\
 \langle\,(\mathbf{C}_{I}+\mathbf{B}^{real}_{n})X^{img},X^{img}\rangle\,+\langle\,\mathbf{B}^{img}_{n}X^{real},X^{img}\rangle&=&\langle\,Y^{img},X^{img}\rangle\label{systemsolve1-small-sub2-ac-2}.
\end{eqnarray}
Summing up \eqref{systemsolve1-small-sub1-ac-2} and \eqref{systemsolve1-small-sub2-ac-2} we obtain
\begin{equation}
\begin{split}
\langle\,\mathbf{C}_{I}X^{real},X^{real}\rangle\,+\langle\,\mathbf{B}^{real}_{n}X^{real},X^{real}\rangle\,+\langle\,\mathbf{C}_{I}X^{img},X^{img}\rangle\,+\langle\,\mathbf{B}^{real}_{n}X^{img},X^{img}\rangle\,\\
=\langle\,Y^{real},X^{real}\rangle+\langle\,Y^{img},X^{img}\rangle.\label{systemsolve1-small-sub1-2-ac-3}
\end{split}
\end{equation}
The right-hand side in \eqref{systemsolve1-small-sub1-2-ac-3} can be estimated as
\begin{equation}\label{systemsolve1-small-sub1-2-ac-4}
\begin{split}
\langle\,X^{real},X^{real}\rangle^{1\slash\,2}\langle\,Y^{real},Y^{real}\rangle^{1\slash\,2}+\langle\,X^{img},X^{img}\rangle^{1\slash\,2}\langle\,Y^{img},Y^{img}\rangle^{1\slash\,2}\\
\leq2\langle\,X^{|\cdot|},X^{|\cdot|}\rangle^{1\slash\,2}\langle\,Y^{|\cdot|},Y^{|\cdot|}\rangle^{1\slash\,2}.
\end{split}
\end{equation}

Here, $X^{|\cdot|}:=(\|X^{real}\|^2+\|X^{img}\|^2)^\frac{1}{2}$. Let us now consider the right hand side in \eqref{systemsolve1-small-sub1-2-ac-3}. First we have
\begin{equation}\label{B_n-estimation}
\vert \langle\; B_n^{real} X^{real},\; X^{real}\rangle \vert \leq \Vert B_n^{real}\Vert_2 \vert X^{real}\vert^2_2
\end{equation}
where $\Vert B_n^{real}\Vert^2_2:=\sum^{M}_{i,\; j =1}\vert ( B^{real}_n)_{i,j}\vert_2^2$ and $(B^{real}_n)_{i,j}:=\Re\; \Gamma^{\omega}(z_i, z_j)$ if $i\ne j$ 
and $(B^{real}_n)_{i,i}:=0$ for $i,j=1,...,M$. 
Hence from (\ref{kupradzeten1}) $\vert (B^{real}_n)_{i,j} \vert_2 \leq \frac{\mathring{c}}{4\pi \vert z_i-z_j\vert },\; i\ne j$. From the observation before Lemma \ref{capacitance-eig-single}, 
we deduce that
\begin{align}\label{estimate-B-n-negative-realpart-1}
\sum^{M}_{i, j=1}(B^{real}_n)^2_{i,j}
&\leq  MM_{max}\left[\frac{{\mathring{c}}^2}{(4\pi)^2 d^2}+\sum^{[\textcolor{black}{\epsilon}^{-\frac{1}{3}}]}_{n=1}
2[(2n+1)^3-(2n-1)^3]\frac{{\mathring{c}}^2}{(4\pi)^2 n^2 \left(\frac{\textcolor{black}{\epsilon}^{\frac{1}{3}}}{2}\right)^{2}}\right]\nonumber\\
&\leq  MM_{max}\left[\frac{{\mathring{c}}^2}{(4\pi)^2  }\epsilon^{-2t}+ 13 {\mathring{c}}^2 \frac{\textcolor{black}{\epsilon}^{-\frac{2}{3}}}{\pi^2}\sum^{[\textcolor{black}{\epsilon}^{-\frac{1}{3}}]}_{n=1}1\right]
=  \frac{ M_{max}^2{\mathring{c}}^2 (13\textcolor{black}{\epsilon}^{-2}+\textcolor{black}{{\epsilon}^{-2t-1}})}{\pi^2}
\end{align}

or
\begin{equation}\label{estimate-B-n-negative-realpart}
\Vert B_n^{real}\Vert_2 \leq \frac{M_{max}\,\mathring{c}}{\pi}(13\textcolor{black}{\epsilon}^{-2}+\textcolor{black}{{\epsilon}^{-2t-1}})^{\frac{1}{2}}
\end{equation}
and then

\begin{equation}\label{B_n-estimation-1}
\vert \langle\; B_n^{real} X^{real},\; X^{real}\rangle \vert +\vert \langle\; B_n^{real} X^{img},\; X^{img}\rangle \vert \leq  \left(\frac{M_{max}\,\mathring{c}}{\pi}(13\textcolor{black}{\epsilon}^{-2}+\textcolor{black}{{\epsilon}^{-2t-1}})^{\frac{1}{2}}\right)\vert X\vert^2_2
\end{equation}

Using Lemma \ref{capacitance-eig-single}, we deduce that
\begin{equation}\label{lower-bound-C-I}
 \langle\,\mathbf{C}_{I}X^{real},X^{real}\rangle\ + \langle\,\mathbf{C}_{I}X^{img},X^{img}\rangle\ \geq (\lambda +2 \mu)^{-1}({\max}C^a_m)^{-1} \vert X \vert^2_2
\end{equation}

From (\ref{systemsolve1-small-sub1-2-ac-3}), (\ref{systemsolve1-small-sub1-2-ac-4}), (\ref{B_n-estimation-1}) and (\ref{lower-bound-C-I}), we deduce that
\begin{equation}\label{invertcondition-quad-approx}
\begin{split}
 \sum_{m=1}^{M}|X_m|^{2}
\leq \left((\lambda +2 \mu)^{-1}({\max}C^a_m)^{-1}-\frac{M_{max}\,\mathring{c}\,(13\textcolor{black}{\epsilon}^{-2}+\textcolor{black}{{\epsilon}^{-2t-1}})^{\frac{1}{2}}}{\pi}\right)^{-2}\sum_{m=1}^{M}\left|Y_m\right|^2
 \end{split}
\end{equation}
if $(\lambda +2 \mu)^{-1}({\max}C^a_m)^{-1}>\frac{M_{max}\,\mathring{c}\,(13\textcolor{black}{\epsilon}^{-2}+\textcolor{black}{{\epsilon}^{-2t-1}})^{\frac{1}{2}}}{\pi}$, which holds true only if $t\leq\frac{1}{2}$.
Hence, \eqref{invertcondition-quad-approx} turns out to be
\begin{equation}\label{invertcondition-quad-approx-1}
\begin{split}
 \sum_{m=1}^{M}|X_m|^{2}
\leq \left((\lambda +2 \mu)^{-1}({\max}C^a_m)^{-1}-\frac{\sqrt{14}M_{max}\,\mathring{c}\,\textcolor{black}{\epsilon}^{-1}}{\pi}\right)^{-2}\sum_{m=1}^{M}\left|Y_m\right|^2
 \end{split}
\end{equation}
if $(\lambda +2 \mu)^{-1}({\max}C^a_m)^{-1}>\frac{\sqrt{14}M_{max}\,\mathring{c}\,\textcolor{black}{\epsilon}^{-1}}{\pi}$ and then the matrix $\mathbf{B}$ 
in algebraic system \eqref{fracqcfracmain} is invertible.

\end{proofe}

\subsection{The limiting model}\label{subsection-piecewise-constant} 
 
From the function $K$, we define a bounded function $K_{\textcolor{black}{\epsilon}}: \mathbb{R}^3\rightarrow \mathbb{R}$ as follows: 
\begin{equation}
K_{\textcolor{black}{\epsilon}}(x):= K_{\textcolor{black}{\epsilon}}(z_m):=
\left\{\begin{array}{ccc}
K(z_m)+1&   \mbox{ if }& x\in \Omega_m \\
0 & \mbox{ if }& x\notin \Omega_m \mbox{ for any } m=1,\dots,[a^{-1}].
\end{array}\right.
\end{equation}
Hence each $\Omega_m$ contains $[ K_{\textcolor{black}{\epsilon}}(z_m)]$ obstacles and $K_{max}:=\sup_{z_m}K_{\textcolor{black}{\epsilon}}(z_m)$.

 \par
Let  ${\bold{C}_{\textcolor{black}{\epsilon}}}$ be the $3 \times 3$ matrix having  entries as piecewise constant functions such that ${\bold{C}_{\textcolor{black}{\epsilon}}}\vert_{\Omega_m}=\bar{C}_m$ for all $m=1,\dots,M$ 
and vanishes outside $\Omega$. Here, $\bar{C}_m$ are the capacitances of $B_m$'s.  From \cite{C-S:2015}, we can observe that $\bar{C}_m$ are defined through
defined through $C_m:=\bar{C}_m \textcolor{black}{\epsilon}$, and are independent of $\textcolor{black}{\epsilon}$.

\par We set 
\begin{equation}\label{sup-cap-all}
\mathcal{C}:=\max\limits_{1\leq{m}\leq{M}}\| \bar{C}_m\|_{\infty}.
\end{equation}
Consider the Lippmann-Schwinger equation
\begin{eqnarray}\label{fracqcfracmain-effect-int}
 Y_{\textcolor{black}{\epsilon}}(z) +\int_{\Omega} \Gamma^{\omega}(z,y)K_{\textcolor{black}{\epsilon}}(y){\bold{C}_{\textcolor{black}{\epsilon}}}(y) Y_{\textcolor{black}{\epsilon}}(y) dy &=&U^{i}(z, \theta), z\in \Omega
\end{eqnarray}
and set the Lam\'e potential
\begin{eqnarray}
 V(Y)(x):=\int_{\Omega}\Gamma^{\omega}(x,y)K_{\textcolor{black}{\epsilon}}(y){\bold{C}_{\textcolor{black}{\epsilon}}}(y)Y(y)dy,\qquad x\in\mathbb{R}^3.
\end{eqnarray}

The coefficients $K_{\textcolor{black}{\epsilon}}$ and $\bold{C}_{\textcolor{black}{\epsilon}}$ are uniformly bounded. The next lemma concerns the mapping properties of the Lam\'e potential. These properties are proved for the scalar Poisson potential in \cite{C-K:1998}, for instance. Similar arguments are applicable for the Lam\'e potential as well, so we omit to give the details.
\begin{lemma}\label{Lip-Schw}
The operator $V:{L}^2(\Omega)\rightarrow {H}^2(\varOmega)$ is well defined and it is a bounded operator for any bounded domain $\varOmega$ in $\mathbb{R}^3$, 
i.e. there exists a positive 
constant $c_0$ such that
\begin{eqnarray}\label{H2normofY}
\|V(Y)\|_{H^{2}(\Omega)}\leq c_0 \|Y\|_{L^{2}(\Omega)}.
\end{eqnarray}

\end{lemma}

\bigskip

 We have also the following lemma.
 \bigskip
 
\begin{lemma}\label{invertibility-of-VC}
 There exists one and only one solution $Y$ of the Lippmann-Schwinger equation (\ref{fracqcfracmain-effect-int}) and it satisfies the estimate
 \begin{eqnarray}\label{est-Lipm-Sch}
  \Vert Y\Vert_{L^\infty(\Omega)}\leq C \Vert U^i\Vert_{H^2(\Omega)}& \mbox{ and }& \Vert \nabla Y\Vert_{L^\infty(\Omega)}\leq C^{\prime} \Vert U^i\Vert_{H^2(\tilde{\Omega})},
 \end{eqnarray}
where $\tilde{\Omega}$ being a large bounded domain which contains $\bar{\Omega}$.
\end{lemma}

\begin{proofe}{\it{Proof of Lemma \ref{invertibility-of-VC}}}

Using the Lemma \ref{Lip-Schw}, we see that $I+V:L^2(\Omega) \longmapsto L^2(\Omega)$ 
is Fredholm with index zero and then we can apply the Fredholm alternative to $I +V:L^2(\Omega) \longmapsto L^2(\Omega)$. 
The uniqueness is a consequence of the uniqueness of the scattering problem corresponding to the model 

\begin{equation}\label{elas-scat-equi}
 (\Delta^e +\omega^2I-K_{\textcolor{black}{\epsilon}} \bold{C}_{\textcolor{black}{\epsilon}}) Y =0,\; \mbox{ in } \mathbb{R}^3
\end{equation}
where $Y:=Y^i +Y^s$ and $Y^s$ satisfies the Kupradze radiation conditions and $Y^i$ is an incident field.

The estimate (\ref{est-Lipm-Sch}) can be derived, as it is done in \cite{DPC-SM13-3} for the acoustic case, by coupling the invertibility of 
$I+V:L^2(\Omega) \longmapsto L^2(\Omega)$ and the $W^{2, p}-$ interior estimates of the solutions of the system $(\Delta^e +\omega^2I-K_{\textcolor{black}{\epsilon}} \bold{C}_{\textcolor{black}{\epsilon}}) Y =0$.

\end{proofe}

\subsubsection{Case when the obstacles are arbitrarily distributed}\label{arbitrarely-distributed}
 
 The capacitances of the obstacles $B_j$, i.e. $C_j$ are bounded by their Lipschitz constants, see \cite{DPC-SM13}, and we assumed that these Lipschitz constants are uniformly bounded.
 Hence $\bold{C}_{\textcolor{black}{\epsilon}}$ is bounded in $L^{2}(\Omega)$ and then there exists a function $\bold{C}_0$ in $L^2(\Omega)$ (actually in every $L^p(\Omega)$) such that $\bold{C}_{\textcolor{black}{\epsilon}}$ converges weakly to $\bold{C}_0$
 in $L^{2}(\Omega)$. Now, since $K$ is continuous hence $K_{\textcolor{black}{\epsilon}}$ converges to $(K+1)$ in $L^\infty(\Omega)$ and hence in $L^2(\Omega)$. 
 Then we can show that $K_{\textcolor{black}{\epsilon}} \bold{C}_{\textcolor{black}{\epsilon}}$ converges to $(K+1) \bold{C}_0$ in $L^2(\Omega)$. 
 \bigskip
 
 Since $K\bold{C}_{\textcolor{black}{\epsilon}}$ is bounded in $L^{\infty}(\Omega)$, then from the invertibility of the Lippmann-Schwinger equation and the mapping properties of the Lam\'e potential, 
 see Lemma \ref{invertibility-of-VC},
 we deduce that $\Vert U_{\textcolor{black}{\epsilon}}^{t}\Vert_{H^2(\Omega)}$ is bounded and in particular, up to a sub-sequence, $U_{\textcolor{black}{\epsilon}}^{t} \,(:=Y_{\textcolor{black}{\epsilon}})$ tends to $U_{0}^{t}$ in $L^2(\Omega)$. 
 From the convergence of $K_{\textcolor{black}{\epsilon}}\bold{C}_{\textcolor{black}{\epsilon}}$ to $(K+1)\bold{C}_0$
 and the one of $U_{\textcolor{black}{\epsilon}}^{t}$ to $U_{0}^{t}$ and (\ref{fracqcfracmain-effect-int}), 
 we derive the following equation satisfied by 
 $U_{0}^{t}(x)$:
 \begin{eqnarray}\label{fracqcfracmain-effect-int0}
 U_{0}^{t}(x) +\int_{\Omega}\Gamma^{\omega}(x, y) (K+1)\bold{C}_0(y) U_{0}^{t}(y)dy=U^i( x, \theta)\; \mbox{ in } \Omega.
 \end{eqnarray}
This is the Lippmann-Schwinger equation corresponding to the scattering problem
\begin{eqnarray}(\Delta^e + \omega^{2}-(K+1) \bold{C}_0)U_{0}^{t}=0 \mbox{ in }\mathbb{R}^{3},\label{piecewise-2}
\end{eqnarray}
with $U_{0}^{t}=U_{0}^s +U^i$,
and $U^s$ satisfies the Kupradze radiation conditions. As the corresponding farfields are of the form
 $$
 U_{0, p}^{\infty}(\hat{x}, \theta)=-\int_{\Omega}\frac{1}{4\pi\,c_p^{2}}(\hat{x}\otimes\hat{x})e^{-i\frac{\omega}{c_p} \hat{x}\cdot y}(K+1)\bold{C}_0(y)U_{0}^{t}(y)dy,
 $$
  $$
 U_{0, s}^{\infty}(\hat{x}, \theta)=-\int_{\Omega}\frac{1}{4\pi\,c_s^{2}}({\rm \textbf{I}}-\hat{x}\otimes\hat{x})e^{-i\frac{\omega}{c_s} \hat{x}\cdot y}(K+1)\bold{C}_0(y)U_{0}^{t}(y)dy
 $$
 and the ones of $U_{\textcolor{black}{\epsilon}}^{t}$ are of the form
 $$
 U_{\textcolor{black}{\epsilon},p}^{\infty}(\hat{x}, \theta)=-\int_{\Omega}\frac{1}{4\pi\,c_p^{2}}(\hat{x}\otimes\hat{x})e^{-i\frac{\omega}{c_p} \hat{x}\cdot y}K_{\textcolor{black}{\epsilon}}\bold{C}_{\textcolor{black}{\epsilon}}(y)U_{\textcolor{black}{\epsilon}}^{t}(y)dy
 $$
  $$
 U^{\infty}_{\textcolor{black}{\epsilon},s}(\hat{x}, \theta)=-\int_{\Omega}\frac{1}{4\pi\,c_s^{2}}({\rm \textbf{I}}-\hat{x}\otimes\hat{x})e^{-i\frac{\omega}{c_s} \hat{x}\cdot y}K_{\textcolor{black}{\epsilon}}\bold{C}_{\textcolor{black}{\epsilon}}(y)U_{\textcolor{black}{\epsilon}}^{t}(y)dy
 $$
 we deduce that
 $$
 U^\infty_{\textcolor{black}{\epsilon},p}(\hat{x},\theta)-U^\infty_{0,p}(\hat{x},\theta)=o(1) \mbox{ and } U^{\infty}_{\textcolor{black}{\epsilon},s}(\hat{x},\theta)-U^\infty_{0,s}(\hat{x},\theta)=o(1),\; {\textcolor{black}{\epsilon}}\rightarrow 0, \mbox{ uniformly in terms of } \hat x , \theta\; \in \mathbb{S}^{2}.
 $$

 \subsubsection{Case when $K$ is H$\ddot{\mbox{o}}$lder continuous}\label{smoothly-distributed}
 
 If we assume that $K\in C^{0, \gamma}(\Omega),\; \gamma \in (0, 1]$, then we have the estimate $\Vert (K+1) -K_{\textcolor{black}{\epsilon}}\Vert_{L^\infty(\Omega)}\leq C {\textcolor{black}{\epsilon}}^{\gamma}$, ${\textcolor{black}{\epsilon}}<<1$.
 Since the capacitances of the obstacles are assumed to be equal, we set $\bold{C}_0$ to be a constant in $\Omega$ and $\bold{C}_0=0$ in $\mathbb{R}^3\setminus \Omega$. 
 Recall that $U_0$ and $U_{\textcolor{black}{\epsilon}}$ are solutions of the Lippmann-Schwinger equations 
 $$
 U_0+\int_{\Omega}(K+1)\Gamma^{\omega}(x, y)\bold{C}_0(y)U_{0}^{t}(y)dy=U^i
 $$
 and 
 $$
 U_{\textcolor{black}{\epsilon}}+\int_{\Omega}K_{\textcolor{black}{\epsilon}}\Gamma^{\omega}(x, y) \bold{C}_0(y)U_{\textcolor{black}{\epsilon}}^{t}(y)dy=U^i.
 $$
 From the estimate $\Vert (K+1) -K_{\textcolor{black}{\epsilon}}\Vert_{L^\infty(\Omega)}\leq C {\textcolor{black}{\epsilon}}^{\gamma}$, ${\textcolor{black}{\epsilon}}<<1$, we derive the estimate
 \begin{equation}\label{appro-0-a}
  U_0^\infty(\hat{x}, \theta)-U^\infty_{\textcolor{black}{\epsilon}}(\hat{x}, \theta)=O({\textcolor{black}{\epsilon}}^\gamma),\; {\textcolor{black}{\epsilon}}<<1, \mbox{ uniformly in terms of } \hat x , \theta\; \in \mathbb{S}^{2}.
 \end{equation}

\subsection{The approximation by the algebraic system}\label{the approximation-by-AS}
For each $m=1,\dots,M$, we rewrite the equation \eqref{fracqcfracmain-effect-int} as follows
\begin{eqnarray}\label{fracqcfracmain-effect-int-1}
 U_{\textcolor{black}{\epsilon}}(z_m) +\sum_{\substack{j=1 \\ j\neq m}}^{M} \Gamma^{\omega}(z_m,z_j)\bar{C}_j U_{\textcolor{black}{\epsilon}}(z_j) {\textcolor{black}{\epsilon}}
 &=&U^{i}(z_m, \theta)+A+B, 
\end{eqnarray}

where 
$$
A:=\sum_{\substack{j=1 \\ j\neq m}}^{[{\textcolor{black}{\epsilon}}^{-1}]} \Gamma^{\omega}(z_m,z_j)K_{\textcolor{black}{\epsilon}}(z_j)\bar{C}_j U_{\textcolor{black}{\epsilon}}(z_j) Vol(\Omega_j)-\int_{\Omega} \Gamma^{\omega}(z_m,y)K_{\textcolor{black}{\epsilon}}(y){\bold{C}_{\textcolor{black}{\epsilon}}}(y) U_{\textcolor{black}{\epsilon}}(y) dy
$$
and 
$$
B:= \sum_{\substack{j=1 \\ j\neq m}}^{M} \Gamma^{\omega}(z_m,z_j)\bar{C}_j U_{\textcolor{black}{\epsilon}}(z_j) {\textcolor{black}{\epsilon}}-\sum_{\substack{j=1 \\ j\neq m}}^{[{\textcolor{black}{\epsilon}}^{-1}]} \Gamma^{\omega}(z_m,z_j)K_{\textcolor{black}{\epsilon}}(z_j)\bar{C}_j U_{\textcolor{black}{\epsilon}}(z_j) Vol(\Omega_j).
$$

Let us estimate the  quantities $A$ and $B$. 
\subsubsection{Estimate of $A$}
By the decomposition of $\Omega$, $\Omega:=\cup^{[{\textcolor{black}{\epsilon}}^{-1}]}_{l=1}$, we have
 \begin{equation}\label{integralonomega}
            \int_{\Omega} \Gamma^{\omega}(z_m,y){K_{\textcolor{black}{\epsilon}}}(y){\bold{C}_{\textcolor{black}{\epsilon}}}(y) U_{\textcolor{black}{\epsilon}}(y) dy=\sum_{l=1}^{[{\textcolor{black}{\epsilon}}^{-1}]}\int_{\Omega_l} \Gamma^{\omega}(z_m,y){K_{\textcolor{black}{\epsilon}}}(y){\bold{C}_{\textcolor{black}{\epsilon}}}(y) U_{\textcolor{black}{\epsilon}}(y) dy.
       \end{equation}
\begin{align}
\mbox{Hence,\quad }A&:= \int_{\Omega_m} \Gamma^{\omega}(z_m,y)K_{\textcolor{black}{\epsilon}}(y){\bold{C}_{\textcolor{black}{\epsilon}}}(y) U_{\textcolor{black}{\epsilon}}(y) dy\nonumber\\
&\quad+\sum_{\substack{j=1 \\ j\neq m}}^{[{\textcolor{black}{\epsilon}}^{-1}]} \left [\Gamma^{\omega}(z_m,z_j)K_{\textcolor{black}{\epsilon}}(z_j)\bar{C}_j U_{\textcolor{black}{\epsilon}}(z_j) Vol(\Omega_j)
                  -\int_{\Omega_j} \Gamma^{\omega}(z_m,y)K_{\textcolor{black}{\epsilon}}(y){\bold{C}_{\textcolor{black}{\epsilon}}}(y) U_{\textcolor{black}{\epsilon}}(y) dy \right].\label{rearranging-A}
\end{align}
 
 For $l\neq m$, we have
 \begin{align}\label{integralonomega-subelements}
  \int_{\Omega_l} \Gamma^{\omega}(z_m,y){K_{\textcolor{black}{\epsilon}}}(y){\bold{C}_{\textcolor{black}{\epsilon}}}(y) U_{\textcolor{black}{\epsilon}}(y) dy &- \Gamma^{\omega}(z_m,z_l){K_{\textcolor{black}{\epsilon}}}(z_l)\bar{C}_l U_{\textcolor{black}{\epsilon}}(z_l) Vol(\Omega_l) \nonumber\\
  &={K_{\textcolor{black}{\epsilon}}}(z_l)\bar{C}_l\int_{\Omega_l} \left[\Gamma^{\omega}(z_m,y) U_{\textcolor{black}{\epsilon}}(y)-\Gamma^{\omega}(z_m,z_l) U_{\textcolor{black}{\epsilon}}(z_l)\right] dy.
 \end{align}
  We set $f(z_m,y)=\Gamma^{\omega}(z_m,y)  U_{\textcolor{black}{\epsilon}}(y)$ then every component $f_i(z_m,y)$ of $f(z_m,y)$ satisfies
  $$f_i(z_m,y)-f_i(z_m,z_l)=(y-z_l)R^i_l(z_m,y)$$
  where
  \begin{eqnarray}\label{taylorremind1}
   R^i_l(z_m,y)
   &=&\int_0^1\nabla_y f_i(z_m,y-\beta(y-z_l))\,d\beta\nonumber\\
   &=&\int_0^1\sum^3_{j=1}\nabla_y \left[\Gamma_{i,j}^{\omega}(z_m,y-\beta(y-z_l)) U_{{\textcolor{black}{\epsilon}},j}(y-\beta(y-z_l))\right]\,d\beta\nonumber\\
   &=&\int_0^1\sum^3_{j=1}\left[\nabla_y\Gamma_{i, j}^{\omega}(z_m,y-\beta(y-z_l))\right] U_{{\textcolor{black}{\epsilon}},j}(y-\beta(y-z_l))\,d\beta\nonumber\\
   &&+\int_0^1\sum^3_{j=1}\Gamma_{i, j}^{\omega}(z_m,y-\beta(y-z_l))\left[\nabla_y U_{{\textcolor{black}{\epsilon}},j}(y-\beta(y-z_l))\right]\,d\beta.
  \end{eqnarray}

 From \eqref{gradkupradzeten1dfnt-1} and from Section \ref{The-relative-distribution-of-the-small-bodies}, we derive for $l \neq m$
   \begin{eqnarray*}
   \vert\Gamma_{i,j}^{\omega}(z_m,y-\beta(y-z_l))\vert\leq\frac{\mathring{c}}{4\pi\; n\frac{{\textcolor{black}{\epsilon}}^{\frac{1}{3}}}{2}},& \mbox{ and }&\vert\nabla_y\Gamma_{i, j}^{\omega}(z_m,y-\beta(y-z_l))\vert\leq\frac{\mathring{c}}{4\pi\;n^2 \left(\frac{{\textcolor{black}{\epsilon}}^{\frac{1}{3}}}{2}\right)^{2}}
  \end{eqnarray*}
  where $\mathring{c}$ depends only on $\omega$ and some universal constants. Then
 \begin{eqnarray}\label{taylorremind-effect}
 \hspace{-.5cm}\vert R_l(z_m,y) \vert
   &\leq & \frac{\mathring{c}}{2\pi\;n {\textcolor{black}{\epsilon}}^{\frac{1}{3}}}\left(\frac{{1}}{n{\textcolor{black}{\epsilon}}^{\frac{1}{3}}}\int_0^1 {\vert U_{\textcolor{black}{\epsilon}}(y-\beta(y-z_l))\vert}d\beta+\int_0^1{\vert\nabla_y U_{\textcolor{black}{\epsilon}}(y-\beta(y-z_l))\vert}d\beta\right).
  \end{eqnarray}

  Then, for $l\neq m$, \eqref{integralonomega-subelements} and \eqref{taylorremind-effect} and observing that $\bar{C}_l$ is a constant matrix in $\Omega_l$, imply the estimate
  
\begin{eqnarray}\label{integralonomega-subelements-abs}
  \Big\vert\int_{\Omega_l}   &    \Gamma^{\omega}(z_m,y) {K_{\textcolor{black}{\epsilon}}}(y){\bold{C}_{\textcolor{black}{\epsilon}}}(y) U_{\textcolor{black}{\epsilon}}(y) dy-\Gamma^{\omega}(z_m,z_l) & {K_{\textcolor{black}{\epsilon}}}(z_l)\bar{C}_l U_{\textcolor{black}{\epsilon}}(z_l) Vol(\Omega_l)\Big\vert 
  \nonumber \\
  \leq\,  &\frac{\mathring{c} \bar{C}_l K_{\textcolor{black}{\epsilon}}(z_l)}{\pi\;n^2 {\textcolor{black}{\epsilon}}^{\frac{2}{3}}}\int_{\Omega_l} \Big[\int_0^1 \vert U_{\textcolor{black}{\epsilon}}  (y-\beta(y-z_l))\vert\;d\beta\Big] &     \vert y-z_l\vert dy
    \nonumber\\
               &&
              + \frac{\mathring{c} \bar{C}_l K_{\textcolor{black}{\epsilon}}(z_l)}{2\pi\;n {\textcolor{black}{\epsilon}}^{\frac{1}{3}}}   \int_{\Omega_l} \left[\int_0^1 {\vert \nabla_y U_{\textcolor{black}{\epsilon}}(y-\beta(y-z_l))\vert}d\beta \right]\vert y-z_l\vert dy 
 \nonumber\\ 
    \substack{\leq\\ \eqref{est-Lipm-Sch}}&c_1\frac{\textcolor{black}{[{K_{\textcolor{black}{\epsilon}}}(z_l)]}\bar{C}_l}{n^2 {\textcolor{black}{\epsilon}}^{\frac{2}{3}}}\,{\textcolor{black}{\epsilon}}^\frac{4}{3}
                    \,\substack{\leq\\ \eqref{sup-cap-all}}\, c_1\frac{\textcolor{black}{K_{\max}}\mathcal{C}}{n^2 }\,{\textcolor{black}{\epsilon}}^\frac{2}{3},\qquad
 \end{eqnarray} 
 \noindent
for a suitable constant $c_1$.
 \par Regarding the integral $\int_{\Omega_m} \Gamma^{\omega}(z_m,y){\bold{C}_{\textcolor{black}{\epsilon}}}(y) U_{\textcolor{black}{\epsilon}}(y) dy$ we do the following estimates:
\begin{eqnarray}\label{estmatemthint-effe-acc}
 \Big\vert\int_{\Omega_m}& \Gamma^{\omega}(z_m,y){K_{\textcolor{black}{\epsilon}}}(y)  {\bold{C}_{\textcolor{black}{\epsilon}}}(y) U_{\textcolor{black}{\epsilon}}(y) dy\Big\vert_{\infty} \nonumber\\
 \substack{\leq\\ \eqref{est-Lipm-Sch}}&c_2{K_{\textcolor{black}{\epsilon}}}(z_m)\bar{C}_m\int_{\Omega_m} \vert \Gamma^{\omega}(z_m,y)\vert_{\infty} dy \nonumber\\
 \substack{\leq\\ \eqref{gradkupradzeten1dfnt-1}} &\frac{\mathring{c}}{4\pi} c_2 {K_{\textcolor{black}{\epsilon}}}(z_m)\bar{C}_m\Big(\int_{B(z_m,r)} \frac{1}{|z_m-y|} &\hspace{-0.5cm}dy   +\int_{\Omega_m\setminus B(z_m,r)} \frac{1}{|z_m-y|} dy\Big)\qquad\nonumber\\
&&{\Big(\mbox{here, }\frac{1}{|z_m-y|} \in L^1(B(z_m,r)), r<\frac{{\textcolor{black}{\epsilon}}^{\frac{1}{3}}}{2} \Big)}\nonumber\\
 \leq\,& \frac{\mathring{c}c_2}{4\pi}{K_{\textcolor{black}{\epsilon}}}(z_m)\bar{C}_m\Big(\sigma(\mathbb{S}^{3-1})\int_0^r \frac{1}{s}s^{3-1}&\hspace{-0.5cm} ds  +\frac{1}{r}Vol(\Omega_m\setminus B(z_m,r))\Big)\nonumber\\
 =\,& \underbrace{\left(2\pi r^2 +\frac{1}{r}\left[{\textcolor{black}{\epsilon}}-\frac{4}{3}\pi r^3\right]\right)}_{=:lm(r,{\textcolor{black}{\epsilon}})}\frac{\mathring{c}}{4\pi} c_2\bar{C}_m{K_{\textcolor{black}{\epsilon}}}&\hspace{-0.3cm}(z_m)\nonumber\\
 \leq\,&\frac{\mathring{c}}{4\pi}c_2{K_{\textcolor{black}{\epsilon}}}(z_m)\bar{C}_m~lm(r^c,{\textcolor{black}{\epsilon}}),\qquad\qquad\nonumber\\
  &&\mbox{here $r^c$ is the value of $r$ where $lm(r,{\textcolor{black}{\epsilon}})$ attains maximum}.\nonumber\\
 &&{\partial_r lm(r,{\textcolor{black}{\epsilon}})=0 \Rightarrow 4\pi r-\frac{{\textcolor{black}{\epsilon}}}{r^2}-\frac{8}{3}\pi r=0\Rightarrow  r_c=\left(\frac{3}{4}\pi {\textcolor{black}{\epsilon}}\right)^\frac{1}{3}}\nonumber\\
 &&{\begin{array}{ccc}
                    lm(r_c,{\textcolor{black}{\epsilon}})&=&2\pi\left(\frac{3}{4}\pi\right)^\frac{2}{3} {\textcolor{black}{\epsilon}}^\frac{2}{3}+\left(\frac{4}{3\pi} \right)^\frac{1}{3}{\textcolor{black}{\epsilon}}^\frac{2}{3}-\frac{4}{3}\pi\left(\frac{3}{4}\pi \right)^\frac{2}{3}{\textcolor{black}{\epsilon}}^\frac{2}{3}\\
                    &&\\
                    &=&\left[\frac{2}{3\pi}\left(\frac{3}{4}\pi\right)^\frac{2}{3}+\left(\frac{4}{3\pi} \right)^\frac{1}{3} \right]{\textcolor{black}{\epsilon}}^\frac{2}{3}=\frac{3}{2}\left(\frac{4}{3\pi} \right)^\frac{1}{3}{\textcolor{black}{\epsilon}}^\frac{2}{3}
                   \end{array}
}\nonumber\\ 
\leq\,&\frac{3}{8\pi}\mathring{c}c_2K_{max}\mathcal{C}\left(\frac{4}{3\pi} \right)^\frac{1}{3}{\textcolor{black}{\epsilon}}^\frac{2}{3}.\qquad\qquad
\end{eqnarray}

From \eqref{rearranging-A}, we can have

\begin{align*}
\vert A\vert_{\infty} \leq&\, \vert \int_{\Omega_m} \Gamma^{\omega}(z_m,y)K_{\textcolor{black}{\epsilon}}(y){\bold{C}_{\textcolor{black}{\epsilon}}}(y) U_{\textcolor{black}{\epsilon}}(y) dy \vert_\infty\\
 &\quad+\sum_{\substack{j=1 \\ j\neq m}}^{[{\textcolor{black}{\epsilon}}^{-1}]} \left[ \vert \Gamma^{\omega}(z_m,z_j)K_{\textcolor{black}{\epsilon}}(z_j)\bar{C}_j U_{\textcolor{black}{\epsilon}}(z_j) Vol(\Omega_j)-\int_{\Omega_j} \Gamma^{\omega}(z_m,y)K_{\textcolor{black}{\epsilon}}(y){\bold{C}_{\textcolor{black}{\epsilon}}}(y) U_{\textcolor{black}{\epsilon}}(y) dy \vert_\infty\right].
\end{align*}
which we can estimate by
\begin{align*}
\vert A\vert_{\infty} \leq 
\sum^{[{\textcolor{black}{\epsilon}}^{-\frac{1}{3}}]}_{n=1}2[(2n&+1)^3-(2n-1)^3]
\bigg[ \vert \Gamma^{\omega}(z_m,z_j)K_{\textcolor{black}{\epsilon}}(z_j)\bar{C}_j U_{\textcolor{black}{\epsilon}}(z_j) Vol(\Omega_j)\\
&-\int_{\Omega_j} \Gamma^{\omega}(z_m,y)K_{\textcolor{black}{\epsilon}}(y){\bold{C}_{\textcolor{black}{\epsilon}}}(y) U_{\textcolor{black}{\epsilon}}(y) dy \vert_\infty\bigg]+\vert \int_{\Omega_m} \Gamma^{\omega}(z_m,y)K_{\textcolor{black}{\epsilon}}(y){\bold{C}_{\textcolor{black}{\epsilon}}}(y) U_{\textcolor{black}{\epsilon}}(y) dy \vert_\infty.
\end{align*}
and then
$$
\vert A\vert_{\infty} \leq c_3\mathcal{C}K_{max}[{\textcolor{black}{\epsilon}}^{\frac{2}{3}}+{\textcolor{black}{\epsilon}}^{\frac{1}{3}}].
$$
Finally 
$$
\vert A \vert \leq c_4\mathcal{C}K_{max} {\textcolor{black}{\epsilon}}^{\frac{1}{3}}.
$$
\subsubsection{Estimate of $B$}
\begin{align*}
\sum_{\substack{j=1 \\ j\neq m}}^{M} \Gamma^{\omega}(z_m,z_j)&\bar{C}_j U_{\textcolor{black}{\epsilon}}(z_j){\textcolor{black}{\epsilon}} -\sum_{\substack{j=1 \\ j\neq m}}^{[{\textcolor{black}{\epsilon}}^{-1}]} \Gamma^{\omega}(z_m,z_j)K_{\textcolor{black}{\epsilon}}(z_j)\bar{C}_j U_{\textcolor{black}{\epsilon}}(z_j) Vol(\Omega_j)           \\
=\sum_{\substack{l=1 \\ l\neq m\\ z_l \in \Omega_m}}^{\textcolor{black}{[K_{\textcolor{black}{\epsilon}}(z_m)]}}\Gamma^{\omega}(z_m&,z_l)\bar{C}_l U_{\textcolor{black}{\epsilon}}(z_l) {\textcolor{black}{\epsilon}}+
\sum_{\substack{j=1 \\ j\neq m}}^{[{\textcolor{black}{\epsilon}}^{-1}]} \sum_{\substack{l=1 \\ z_l \in \Omega_j}}^{\textcolor{black}{[K_{\textcolor{black}{\epsilon}}(z_j)]}}\Gamma^{\omega}(z_m,z_l)\bar{C}_l U_{\textcolor{black}{\epsilon}}(z_l) {\textcolor{black}{\epsilon}}-
\sum_{\substack{j=1 \\ j\neq m}}^{[{\textcolor{black}{\epsilon}}^{-1}]} \Gamma^{\omega}(z_m,z_j)K_{\textcolor{black}{\epsilon}}(z_j)\bar{C}_j U_{\textcolor{black}{\epsilon}}(z_j) Vol(\Omega_j)              \\
=\bar{C}_m {\textcolor{black}{\epsilon}}\sum_{\substack{l=1 \\ l\neq m\\ z_l \in \Omega_m}}^{\textcolor{black}{[K_{\textcolor{black}{\epsilon}}(z_m)]}}\Gamma^{\omega}&(z_m,z_l) U_{\textcolor{black}{\epsilon}}(z_l) +\sum_{\substack{j=1 \\ j\neq m}}^{[{\textcolor{black}{\epsilon}}^{-1}]}\bar{C}_j {\textcolor{black}{\epsilon}}\big[\big(\sum_{\substack{l=1 \\ z_l \in \Omega_j}}^{\textcolor{black}{[K_{\textcolor{black}{\epsilon}}(z_j)]}}\Gamma^{\omega}(z_m,z_l) U_{\textcolor{black}{\epsilon}}(z_l)\big)-
 \Gamma^{\omega}(z_m,z_j)\textcolor{black}{[K_{\textcolor{black}{\epsilon}}(z_j)]} U_{\textcolor{black}{\epsilon}}(z_j)\big],
 \end{align*}
 since $ Vol(\Omega_j)={\textcolor{black}{\epsilon}}\frac{\textcolor{black}{[K_{\textcolor{black}{\epsilon}}(z_j)]}}{K_{\textcolor{black}{\epsilon}}(z_j)}\; \mbox{ and } \bar{C}_l=\bar{C}_j, \mbox{ for } l=1, ...,\; K_{\textcolor{black}{\epsilon}}(z_j)$. We write,

\begin{eqnarray}\label{Ej1}
E^j_1&:=&\sum_{\substack{l=1 \\ l\neq m\\ z_l \in \Omega_m}}^{\textcolor{black}{[K_{\textcolor{black}{\epsilon}}(z_m)]}}\Gamma^{\omega}(z_m,z_l) U_{\textcolor{black}{\epsilon}}(z_l)
\end{eqnarray}
 and
\begin{eqnarray}\label{Ej2}
E^j_2&:=&\big[\big(\sum_{\substack{l=1 \\ z_l \in \Omega_j}}^{\textcolor{black}{[K_{\textcolor{black}{\epsilon}}(z_j)]}}\Gamma^{\omega}(z_m,z_l) U_{\textcolor{black}{\epsilon}}(z_l)\big)-
 \Gamma^{\omega}(z_m,z_j)\textcolor{black}{[K_{\textcolor{black}{\epsilon}}(z_j)]} U_{\textcolor{black}{\epsilon}}(z_j)\big]\nonumber\\
 &=&\sum_{\substack{l=1 \\ z_l \in \Omega_j}}^{\textcolor{black}{[K_{\textcolor{black}{\epsilon}}(z_j)]}}\big(\Gamma^{\omega}(z_m,z_l) U_{\textcolor{black}{\epsilon}}(z_l)-
 \Gamma^{\omega}(z_m,z_j) U_{\textcolor{black}{\epsilon}}(z_j)\big).
\end{eqnarray}

We need to estimate $\bar{C}_m {\textcolor{black}{\epsilon}} E^j_1$ and $\sum_{\substack{j=1 \\ j\neq m}}^{[{\textcolor{black}{\epsilon}}^{-1}]}\bar{C}_j {\textcolor{black}{\epsilon}}E^j_2$. \\

\bigskip
Let us write, $f'(z_m,y):=\Gamma^{\omega}(z_m,y) U_{\textcolor{black}{\epsilon}}(y)$. For $z_l\in\Omega_j,\,j\neq m$, using Taylor series, we can write
  $$f'(z_m,z_j)-f'(z_m,z_l)=(z_j-z_l)R'(z_m;z_j,z_l),$$
  with
 \begin{eqnarray}\label{taylorremind1'}
   R'(z_m;z_j,z_l)
   &=&\int_0^1\nabla_y f'(z_m,z_j-\beta(z_j-z_l))\,d\beta.
  \end{eqnarray}
By doing the computations similar to the ones we have performed in (\ref{taylorremind1}-\ref{taylorremind-effect}) and by using Lemma \ref{invertibility-of-VC}, we obtain
\begin{eqnarray}
  \vert \sum_{\substack{j=1 \\ j\neq m}}^{[{\textcolor{black}{\epsilon}}^{-1}]}\bar{C}_j {\textcolor{black}{\epsilon}} E^j_2 \vert \leq c_4\mathcal{C}K_{max} {\textcolor{black}{\epsilon}}^{\frac{1}{3}}   \label{este2j1}
\end{eqnarray}
One can easily see that,
\begin{eqnarray}
\vert\bar{C}_m {\textcolor{black}{\epsilon}} E^j_1\vert \leq \frac{\mathring{c}c_2(K_{max}-1)\mathcal{C}}{4\pi}\frac{{\textcolor{black}{\epsilon}}}{d} = \frac{\mathring{c}c_2(K_{max}-1)\mathcal{C}}{4\pi}
{\textcolor{black}{\epsilon}}^{1-t}.   \label{este1j}
\end{eqnarray}
\bigskip
\subsubsection{End of the approximation by the algebraic system}\label{the approximation-by-AS-SIUBSUB}

Substitution of \eqref{integralonomega} in \eqref{fracqcfracmain-effect-int-1} and using the estimates 
\eqref{integralonomega-subelements-abs} and \eqref{estmatemthint-effe-acc}  associated to $A$ and the estimates 
\eqref{este2j1} and \eqref{este1j} associated to $B$  gives us

\begin{eqnarray}
 U_{\textcolor{black}{\epsilon}}(z_m) +\sum_{\substack{j=1 \\ j\neq m}}^{M} \Gamma^{\omega}(z_m,z_j)\bar{C}_j U_{\textcolor{black}{\epsilon}}(z_j) {\textcolor{black}{\epsilon}}
  &=&U^{i}(z_m, \theta)+O\left(c_4 K_{max} {\textcolor{black}{\epsilon}}^{\frac{1}{3}}\right)+O\left(\frac{\mathring{c}c_2(K_{max}-1)\mathcal{C} }{4\pi} {\textcolor{black}{\epsilon}}^{1-t}\right).\quad\label{fracqcfracmain-effect-int-3}
\end{eqnarray}

We rewrite the algebraic system (\ref{fracqcfracmain}) as 
\begin{equation}\label{alge-rewriten-U}
U_{{\textcolor{black}{\epsilon}},m} + \sum_{\substack{j=1 \\ j\neq m}}^{M} \Gamma^{\omega}(z_m,z_j)\bar{C}_j U_{{\textcolor{black}{\epsilon}}, j} {\textcolor{black}{\epsilon}}
  =U^{i}(z_m)
\end{equation} 
where we set $U_{\textcolor{black}{\epsilon},m}:=-C^{-1}_m Q_m$, recalling that $C_m=\bar{C}_m\; {\textcolor{black}{\epsilon}}$.
\bigskip

Taking the difference between \eqref{fracqcfracmain-effect-int-3} and 
\eqref{alge-rewriten-U}  produces the algebraic system
\begin{eqnarray}\label{fracqcfracmain-effect-int-4}
 (U_{\textcolor{black}{\epsilon},m}-U_{\textcolor{black}{\epsilon}}(z_m)) +\sum_{\substack{j=1 \\ j\neq m}}^{M} \Gamma^{\omega}(z_m,z_j)\bar{C}_j (U_{{\textcolor{black}{\epsilon}},j}-U_{\textcolor{black}{\epsilon}}(z_j)) {\textcolor{black}{\epsilon}} 
 &=&O\left(\mathcal{C} K_{max}({\textcolor{black}{\epsilon}}^{\frac{1}{3}}+{\textcolor{black}{\epsilon}}^{1- t}) \right).\nonumber
\end{eqnarray}

Comparing this system with \eqref{fracqcfracmain} and by using Lemma \ref{Mazyawrkthmela}, we obtain the estimate

\begin{eqnarray}\label{mazya-fnlinvert-small-ac-3-effect-dif}
 \sum_{m=1}^{M}(U_{\textcolor{black}{\epsilon},m}-U_{\textcolor{black}{\epsilon}}(z_m))&=&O\left(\mathcal{C} K_{max} M ({\textcolor{black}{\epsilon}}^{\frac{1}{3}}+{\textcolor{black}{\epsilon}}^{1- t}) \right).
\end{eqnarray}

For the special case $d={\textcolor{black}{\epsilon}}^t,\,M=O({\textcolor{black}{\epsilon}}^{-1})$ with  $t>0$, we have the following approximation of the far-field from the Foldy-Lax asymptotic expansion \eqref{x oustdie1 D_m farmainp-near} and from the definitions $U_{\textcolor{black}{\epsilon},m}:=-C^{-1}_m Q_m$ and $C_m:=\bar{C}_m {\textcolor{black}{\epsilon}}$, for $m=1,\dots,M$:
\begin{eqnarray}\label{x oustdie1 D_m farmain-recent**-effect}
4 \pi c_p^2U_p^\infty(\hat{x},\theta)\cdot \hat{x} &=&-\sum_{j=1}^{M}e^{-i\kappa\hat{x}\cdot z_j}\bar{C}_jU_{{\textcolor{black}{\epsilon}},j}\cdot \hat{x}\; {\textcolor{black}{\epsilon}}\\ \nonumber
&&+O\left({\textcolor{black}{\epsilon}}\hspace{-.03cm}+\hspace{-.03cm}{\textcolor{black}{\epsilon}}^{2-5t+3t\alpha}\hspace{-.03cm}+\hspace{-.03cm}{\textcolor{black}{\epsilon}}^{3-9t+6t\alpha}
\hspace{-.03cm}+\hspace{-.03cm}{\textcolor{black}{\epsilon}}^{1-2t\alpha}\hspace{-.03cm}+
\hspace{-.03cm}{\textcolor{black}{\epsilon}}^{1-3t\alpha}\hspace{-.03cm}+
\hspace{-.03cm}{\textcolor{black}{\epsilon}}^{2-5t+2t\alpha}\right).
\end{eqnarray}
Consider the far-field of type:
\begin{eqnarray}\label{acoustic-farfield-effect}
U^\infty_{\bold{C}_{\textcolor{black}{\epsilon}}}(\hat{x},\theta) &=&- \frac{1}{4\pi c^2_p}(\hat{x} \otimes \hat{x})\int_{\Omega} e^{-i\frac{\omega}{c_p}\hat{x}\cdot{y}}{K_{\textcolor{black}{\epsilon}}}(y){\bold{C}_{\textcolor{black}{\epsilon}}}(y) U_{\textcolor{black}{\epsilon}}(y) dy -\frac{1}{4\pi c^2_s}(I-\hat{x} \otimes \hat{x})\int_{\Omega} e^{-i\frac{\omega}{c_s}\hat{x}\cdot{y}}{K_{\textcolor{black}{\epsilon}}}(y){\bold{C}_{\textcolor{black}{\epsilon}}}(y) U_{\textcolor{black}{\epsilon}}(y) dy. \nonumber
\end{eqnarray}
corresponding to the scattering problem (\ref{elas-scat-equi}) and set 
\begin{equation}\label{u^{infty}_p}
U^\infty_{\bold{C}_{\textcolor{black}{\epsilon}},p}(\hat{x},\theta):= -\frac{1}{4\pi c_p^2}(\hat{x} \otimes \hat{x})\int_{\Omega} e^{-i\frac{\omega}{c_p}\hat{x}\cdot{y}}{K_{\textcolor{black}{\epsilon}}}(y){\bold{C}_{\textcolor{black}{\epsilon}}}(y) U_{\textcolor{black}{\epsilon}}(y) dy 
\end{equation}
and
\begin{equation}\label{u^{infty}_s}
U^\infty_{\bold{C}_{\textcolor{black}{\epsilon}},s}(\hat{x},\theta) := -\frac{1}{4\pi c_s^2}(I- \hat{x} \otimes \hat{x})\int_{\Omega} e^{-i\frac{\omega}{c_s}\hat{x}\cdot{y}}{K_{\textcolor{black}{\epsilon}}}(y){\bold{C}_{\textcolor{black}{\epsilon}}}(y) U_{\textcolor{black}{\epsilon}}(y) dy 
\end{equation}

Taking the difference between  \eqref{x oustdie1 D_m farmain-recent**-effect} and \eqref{acoustic-farfield-effect}  we have:
\begin{eqnarray}\label{acoustic-difference-farfield-effect}
 &&4\pi c_p^2(U^\infty_p(\hat{x},\theta)-U_{\bold{C}_{\textcolor{black}{\epsilon}},p}^\infty(\hat{x},\theta))\cdot \hat{x} \nonumber\\
 &=& \int_{\Omega} e^{-i\kappa_p\hat{x}\cdot{y}}{K_{\textcolor{black}{\epsilon}}}(y){\bold{C}_{\textcolor{black}{\epsilon}}}(y) U_{\textcolor{black}{\epsilon}}(y)\cdot \hat{x} dy- \sum_{j=1}^{M}e^{-i\kappa_p\hat{x}\cdot z_j}\bar{C}_jU_{{\textcolor{black}{\epsilon}},j}\cdot \hat{x}a\nonumber\\ \nonumber
 &&\hspace{0cm}+O\left({\textcolor{black}{\epsilon}}\hspace{-.03cm}+\hspace{-.03cm}{\textcolor{black}{\epsilon}}^{2-5t+3t\alpha}\hspace{-.03cm}+
 \hspace{-.03cm}{\textcolor{black}{\epsilon}}^{3-9t+6t\alpha}\hspace{-.03cm}+\hspace{-.03cm}{\textcolor{black}{\epsilon}}^{1-2t\alpha}\hspace{-.03cm}+
 \hspace{-.03cm}{\textcolor{black}{\epsilon}}^{1-3t\alpha}\hspace{-.03cm}+\hspace{-.03cm}{\textcolor{black}{\epsilon}}^{2-5t+2t\alpha}\right)
\\ \nonumber &=& \sum_{j=1}^{[{\textcolor{black}{\epsilon}}^{-1}]}\int_{\Omega_j} e^{-i\frac{\omega}{c_p}\hat{x}\cdot{y}}{K_{\textcolor{black}{\epsilon}}}(y){\bold{C}_{\textcolor{black}{\epsilon}}}(y) U_{\textcolor{black}{\epsilon}}(y)\cdot \hat{x}dy -\sum_{j=1}^{[{\textcolor{black}{\epsilon}}^{-1}]} \sum_{\substack{l=1 \\ z_l \in \Omega_j}}^{\textcolor{black}{[K_{\textcolor{black}{\epsilon}}(z_j)]}}e^{-i\frac{\omega}{c_p}\hat{x}\cdot z_l}\bar{C}_lU_{{\textcolor{black}{\epsilon}},l}\cdot \hat{x}a\nonumber\\ \nonumber
 &&\hspace{0cm}+O\left({\textcolor{black}{\epsilon}} \hspace{-.03cm}+\hspace{-.03cm}{\textcolor{black}{\epsilon}}^{2-5t+3t\alpha}\hspace{-.03cm}+
 \hspace{-.03cm}{\textcolor{black}{\epsilon}}^{3-9t+6t\alpha}\hspace{-.03cm}+\hspace{-.03cm}{\textcolor{black}{\epsilon}}^{1-2t\alpha}\hspace{-.03cm}+\hspace{-.03cm}{\textcolor{black}{\epsilon}}^{1-3t\alpha}\hspace{-.03cm}+
 \hspace{-.03cm}{\textcolor{black}{\epsilon}}^{2-5t+2t\alpha}\right)
\\ \nonumber &=& \sum_{j=1}^{[{\textcolor{black}{\epsilon}}^{-1}]}{K_{\textcolor{black}{\epsilon}}}(z_j)\bar{C}_j\int_{\Omega_j} \left[e^{-i\frac{\omega}{c_p}\hat{x}\cdot{y}} U_{\textcolor{black}{\epsilon}}(y)\cdot \hat{x} - e^{-i\frac{\omega}{c_p}\hat{x}\cdot z_j}U_{\textcolor{black}{\epsilon}}(z_j)\cdot \hat{x}\right]dy\nonumber\\ \nonumber
\\ \nonumber &&+ \sum_{j=1}^{[{\textcolor{black}{\epsilon}}^{-1}]}\bar{C}_j{\textcolor{black}{\epsilon}} \left[\sum_{\substack{l=1 \\ z_l \in \Omega_j}}^{\textcolor{black}{[K_{\textcolor{black}{\epsilon}}(z_j)]}}\left(e^{-i\frac{\omega}{c_p}\hat{x}\cdot z_j} U_{\textcolor{black}{\epsilon}}(z_j)\cdot \hat{x}-e^{-i\frac{\omega}{c_p}\hat{x}\cdot z_l}U_{\textcolor{black}{\epsilon}}(z_l)\right)+\sum_{\substack{l=1 \\ z_l \in \Omega_j}}^{\textcolor{black}{[K_{\textcolor{black}{\epsilon}}(z_j)]}} e^{-i\frac{\omega}{c_p}\hat{x}\cdot z_l} \left(U_{\textcolor{black}{\epsilon}}(z_l)-U_{{\textcolor{black}{\epsilon}},l}\right)\cdot \hat{x}\right]\nonumber\\ \nonumber
 &&\hspace{0cm}+O\left({\textcolor{black}{\epsilon}}\hspace{-.03cm}+\hspace{-.03cm}{\textcolor{black}{\epsilon}}^{2-5t+3t\alpha}\hspace{-.03cm}+\hspace{-.03cm}
 {\textcolor{black}{\epsilon}}^{3-9t+6t\alpha}\hspace{-.03cm}+\hspace{-.03cm}{\textcolor{black}{\epsilon}}^{1-2t\alpha}
 \hspace{-.03cm}+\hspace{-.03cm}{\textcolor{black}{\epsilon}}^{1-3t\alpha}\hspace{-.03cm}+
 \hspace{-.03cm}{\textcolor{black}{\epsilon}}^{2-5t+2t\alpha}\right)
 \\ \nonumber &=& \sum_{j=1}^{[{\textcolor{black}{\epsilon}}^{-1}]}\int_{\Omega_j}{K_{\textcolor{black}{\epsilon}}}(z_j)\bar{C}_j \left[e^{-i\frac{\omega}{c_p}\hat{x}\cdot{y}} U_{\textcolor{black}{\epsilon}}(y)\cdot \hat{x} - e^{-i\frac{\omega}{c_p}\hat{x}\cdot z_j}U_{\textcolor{black}{\epsilon}}(z_j)\cdot \hat{x}\right]dy\nonumber\\
 &&+ \sum_{j=1}^{[{\textcolor{black}{\epsilon}}^{-1}]}\bar{C}_j{\textcolor{black}{\epsilon}} \sum_{\substack{l=1 \\ z_l \in \Omega_j}}^{\textcolor{black}{[K_{\textcolor{black}{\epsilon}}(z_j)]}}\left(e^{-i\frac{\omega}{c_p}\hat{x}\cdot z_j} U_{\textcolor{black}{\epsilon}}(z_j)\cdot \hat{x}-e^{-i\frac{\omega}{c_p}\hat{x}\cdot z_l}U_{\textcolor{black}{\epsilon}}(z_l)\cdot \hat{x}\right)+\sum_{j=1}^{M}e^{-i\frac{\omega}{c_p}\hat{x}\cdot z_j}\bar{C}_j{\textcolor{black}{\epsilon}} \left[U_{\textcolor{black}{\epsilon}}(z_j)-U_{{\textcolor{black}{\epsilon}},j}\right]\cdot \hat{x}\nonumber\\ \nonumber
 &&\hspace{0cm}+O\left({\textcolor{black}{\epsilon}}\hspace{-.03cm}+\hspace{-.03cm}
 {\textcolor{black}{\epsilon}}^{2-5t+3t\alpha}\hspace{-.03cm}+\hspace{-.03cm}{\textcolor{black}{\epsilon}}^{3-9t+6t\alpha}
 \hspace{-.03cm}+\hspace{-.03cm}{\textcolor{black}{\epsilon}}^{1-2t\alpha}\hspace{-.03cm}+
 \hspace{-.03cm}{\textcolor{black}{\epsilon}}^{1-3t\alpha}\hspace{-.03cm}+\hspace{-.03cm}
 {\textcolor{black}{\epsilon}}^{2-5t+2t\alpha}\right)\\
 &\substack{= \\ \eqref{mazya-fnlinvert-small-ac-3-effect-dif} }& \sum_{j=1}^{[{\textcolor{black}{\epsilon}}^{-1}]}{K_{\textcolor{black}{\epsilon}}}(z_j)\bar{C}_j\int_{\Omega_j} \left[e^{-i\frac{\omega}{c_p}\hat{x}\cdot{y}} U_{\textcolor{black}{\epsilon}}(y)\cdot \hat{x} - e^{-i\frac{\omega}{c_p}\hat{x}\cdot z_j}U_{\textcolor{black}{\epsilon}}(z_j)\cdot \hat{x}\right]dy\nonumber\\
 &&+ \sum_{j=1}^{[{\textcolor{black}{\epsilon}}^{-1}]}\bar{C}_j{\textcolor{black}{\epsilon}} \sum_{\substack{l=1 \\ z_l \in \Omega_j}}^{\textcolor{black}{[K_{\textcolor{black}{\epsilon}}(z_j)]}}\left(e^{-i\frac{\omega}{c_p}\hat{x}\cdot z_j} U_{\textcolor{black}{\epsilon}}(z_j)\cdot \hat{x}-e^{-i\frac{\omega}{c_p}\hat{x}\cdot z_l}U_{\textcolor{black}{\epsilon}}(z_l)\cdot \hat{x}\right)+O\left(\mathcal{C}^2 K_{max}  ({\textcolor{black}{\epsilon}}^{\frac{1}{3}}+{\textcolor{black}{\epsilon}}^{1-t})\right)\nonumber\\
&&+O\left({\textcolor{black}{\epsilon}}\hspace{-.03cm}+\hspace{-.03cm}{\textcolor{black}{\epsilon}}^{2-5t+3t\alpha}\hspace{-.03cm}+
   \hspace{-.03cm}{\textcolor{black}{\epsilon}}^{3-9t+6t\alpha}\hspace{-.03cm}+\hspace{-.03cm}
   {\textcolor{black}{\epsilon}}^{1-2t\alpha}\hspace{-.03cm}+\hspace{-.03cm}{\textcolor{black}{\epsilon}}^{1-3t\alpha}
   \hspace{-.03cm}+\hspace{-.03cm}{\textcolor{black}{\epsilon}}^{2-5t+2t\alpha}\right).
 \end{eqnarray}
Now, let us estimate the  difference $\sum_{j=1}^{[{\textcolor{black}{\epsilon}}^{-1}]}{K_{\textcolor{black}{\epsilon}}}(z_j)\bar{C}_j\int_{\Omega_j} \left[e^{-i\frac{\omega}{c_p}\hat{x}\cdot{y}} U_{\textcolor{black}{\epsilon}}(y) - 
e^{-i\frac{\omega}{c_p}\hat{x}\cdot z_j}U_{\textcolor{black}{\epsilon}}(z_j)\right]dy$.
Write, $f_1(y)=e^{-i\frac{\omega}{c_p}\hat{x}\cdot{y}} U_{\textcolor{black}{\epsilon}}(y)$. Using Taylor series, we can write
  $$f_1(y)-f_1(z_j)=(y-z_j)R_j(y),$$
  with 
  \begin{eqnarray}
   (R_j(y))_{k,l}
   &=&\int_0^1\nabla_y (f_1)_{k,l}(y-\beta(y-z_j))\,d\beta\nonumber\\
   &=&\int_0^1\left[\nabla_y \left[e^{-i\frac{\omega}{c_p}\hat{x}\cdot(y-\beta(y-z_j))} U_{\textcolor{black}{\epsilon}}(y-\beta(y-z_j))\right] \right]_k\,d\beta\nonumber\\
   &=&\int_0^1\left[\nabla_ye^{-i\frac{\omega}{c_p}\hat{x}\cdot(y-\beta(y-z_j))}\right]_k U_{\textcolor{black}{\epsilon}}(y-\beta(y-z_j))\,d\beta\nonumber\\
   &&+\int_0^1e^{-i\frac{\omega}{c_p}\hat{x}\cdot(y-\beta(y-z_j))}\left[\nabla_y U_{\textcolor{black}{\epsilon}}(y-\beta(y-z_j))\right]_k\,d\beta.
  \end{eqnarray}

  We have $\nabla_ye^{-i\frac{\omega}{c_p}\hat{x}\cdot{y}}=-i\frac{\omega}{c_p}\hat{x}e^{-i\frac{\omega}{c_p}\hat{x}\cdot{y}}$ then
  
 \begin{eqnarray}\label{taylorremind-effect-singvariable}
 \vert R_j(y) \vert_\infty
   &\leq&\left(\frac{\omega}{c_p} \int_0^1 |U_{\textcolor{black}{\epsilon}}(y-\beta(y-z_j))|_\infty\, d\beta\,+\,\int_0^1\vert\nabla_y U_{\textcolor{black}{\epsilon}}(y-\beta(y-z_j))\vert_\infty \,d\beta\right).
  \end{eqnarray}
  Using  \eqref{taylorremind-effect-singvariable} we get the estimate
 \begin{equation}\label{integralonomega-subelements-abs-effect}
  \left\vert\sum_{j=1}^{[{\textcolor{black}{\epsilon}}^{-1}]}{K_{\textcolor{black}{\epsilon}}}(z_j)\bar{C}_j\int_{\Omega_j} \left[e^{-i\frac{\omega}{c_p}\hat{x}\cdot{y}}(y) U_{\textcolor{black}{\epsilon}}(y) - e^{-i\frac{\omega}{c_p}\hat{x}\cdot z_j}U_{\textcolor{black}{\epsilon}}(z_j)\right]dy\right \vert_\infty \leq 
\end{equation}

$$
\sum_{j=1}^{[{\textcolor{black}{\epsilon}}^{-1}]}{K_{\textcolor{black}{\epsilon}}}(z_j)\bar{C}_j \left(\frac{\omega}{c_p} \int_{\Omega_j}\vert y-z_j\vert_\infty\int_0^1 |U_{\textcolor{black}{\epsilon}}(y-\beta(y-z_j))|\, d\beta\,dy\right)\, 
$$

$$
  +\,\sum_{j=1}^{[{\textcolor{black}{\epsilon}}^{-1}]}{K_{\textcolor{black}{\epsilon}}}(z_j)\bar{C}_j \left(\int_{\Omega_j}\vert y-z_j\vert\int_0^1\vert\nabla_y U_{\textcolor{black}{\epsilon}}(y-\beta(y-z_j))\vert_\infty \,d\beta\,dy\right)
  $$
  
  $$ 
  \leq \sum_{j=1}^{[{\textcolor{black}{\epsilon}}^{-1}]}{K_{\textcolor{black}{\epsilon}}}(z_j)\bar{C}_jc_1\, {\textcolor{black}{\epsilon}}\,{\textcolor{black}{\epsilon}}^\frac{1}{3}\,\left(\frac{\omega}{c_p} +c_5\right) 
   \leq K_{max} \mathcal{C}c_1\left(\frac{\omega}{c_p} +c_5\right)\,{\textcolor{black}{\epsilon}}^{\frac{1}{3}}.   
$$

 In the similar way, using (\ref{mazya-fnlinvert-small-ac-3-effect-dif}),  we have,
 \begin{eqnarray}\label{integralonomega-subelements-abs-effect-1}
 \left\vert \sum_{j=1}^{[{\textcolor{black}{\epsilon}}^{-1}]}\bar{C}_j{\textcolor{black}{\epsilon}} \sum_{\substack{l=1 \\ z_l \in \Omega_j}}^{\textcolor{black}{[K_{\textcolor{black}{\epsilon}}(z_j)]}}\left(e^{-i\frac{\omega}{c_p}\hat{x}\cdot z_j} U_{\textcolor{black}{\epsilon}}(z_j)-e^{-i\frac{\omega}{c_p}\hat{x}\cdot z_l}U_{\textcolor{black}{\epsilon}}(z_l)\right)\right\vert_\infty
 &\leq& O\left(K_{max} \mathcal{C}({\textcolor{black}{\epsilon}}^{\frac{1}{3}}+{\textcolor{black}{\epsilon}}^{1-t})\right).
 \end{eqnarray}
 
 Using the estimates \eqref{integralonomega-subelements-abs-effect} and \eqref{integralonomega-subelements-abs-effect-1} in \eqref{acoustic-difference-farfield-effect}, we obtain
 
 \begin{eqnarray}\label{acoustic-difference-farfield-effect-1}
{4\pi c^2_p}(U^\infty_p(\hat{x},\theta)-U_{\bold{C}_{\textcolor{black}{\epsilon}},p}^\infty(\hat{x},\theta))\cdot \hat{x} 
 &=& O\left(K_{max}{\textcolor{black}{\epsilon}}^\frac{1}{3} \mathcal{C}c_1\left(\frac{\omega}{c_p} +c_5\right)\right) + O(\mathcal{C}(\mathcal{C}+1)\mathcal{C}K_{max} ({\textcolor{black}{\epsilon}}^{\frac{1}{3}}+{\textcolor{black}{\epsilon}}^{1-t}))\nonumber\\
   &&\hspace{0cm}+O\left({\textcolor{black}{\epsilon}}\hspace{-.03cm}+\hspace{-.03cm}{\textcolor{black}{\epsilon}}^{2-5t+3t\alpha}\hspace{-.03cm}+
   \hspace{-.03cm}{\textcolor{black}{\epsilon}}^{3-9t+6t\alpha}\hspace{-.03cm}+
   \hspace{-.03cm}{\textcolor{black}{\epsilon}}^{1-2t\alpha}\hspace{-.03cm}+
   \hspace{-.03cm}{\textcolor{black}{\epsilon}}^{1-3t\alpha}\hspace{-.03cm}+
   \hspace{-.03cm}{\textcolor{black}{\epsilon}}^{2-5t+2t\alpha}\right) 
 \nonumber\\
    &=&O\left(  {\textcolor{black}{\epsilon}}^{\frac{1}{3}}+{\textcolor{black}{\epsilon}}^{1-t}+{\textcolor{black}{\epsilon}}^{3- 9t+6\alpha t}\hspace{-.03cm}+\hspace{-.03cm}
   {\textcolor{black}{\epsilon}}^{1-3\alpha t}+{\textcolor{black}{\epsilon}}^{2-5t+2t\alpha}\right).
  \end{eqnarray}
 Since $Vol(\Omega)$ is of order ${\textcolor{black}{\epsilon}}^{-1}\Big(\frac{\textcolor{black}{\epsilon}}{2}{{\textcolor{black}{\max\limits_{1\leq m\leq M } \diam (B_m)}}}+\frac{d}{2}\Big)^3$, and $d$ is of the order ${\textcolor{black}{\epsilon}}^t$, we should have {$t\geq\frac{1}{3}$}. Hence
\begin{eqnarray*}
  t\geq\frac{1}{3}; \qquad 1-t>0;   \qquad
      3-9t+6t\alpha>0;\qquad  1-3t\alpha>0;\qquad 2-5t+2t\alpha>0;
  \end{eqnarray*}  
Hence for {$\frac{1}{3}\leq{t}<1$}, we have 
{
\begin{eqnarray*}
  t\alpha<\frac{1}{3};\qquad 2-5t+2t\alpha>0;\qquad
 3-9t+6t\alpha>0.
\end{eqnarray*}
}

Equating $2-5t+2t\alpha=3-9t+6t\alpha$, we find that $\alpha t =t-\frac{1}{4}$ and then $2-5t+2t\alpha=3-9t+6t\alpha=\frac{3}{2}-3t$ and $1-3 \alpha t=\frac{7}{4}-3 t$. 
 In addition, since $\alpha t < \frac{1}{3}$, then $t< \frac{7}{12}$. To make sure that $\frac{3}{2}-3t$ is positive we should have $t<\frac{1}{2}$. Hence, the error is

\begin{eqnarray}\label{elastic-difference-farfield-effect-1**1-1}
{4\pi c^2_p} \left[ U^\infty_p(\hat{x},\theta)-U_{\bold{C}_{\textcolor{black}{\epsilon}},p}^\infty(\hat{x},\theta)\right]\cdot \hat{x}
      &=&O\left( {\textcolor{black}{\epsilon}}^{\frac{1}{3}} +{\textcolor{black}{\epsilon}}^{\frac{3}{2}-3t}\right),\;~~ \frac{1}{3} \leq t < \frac{1}{2}.
  \end{eqnarray}
 
Doing similar computations for the $s$-parts of the far fields as it was done for the $p$-parts, between (\ref{x oustdie1 D_m farmain-recent**-effect}-- \ref{acoustic-difference-farfield-effect-1}), we obtain the following estimates
\begin{eqnarray}\label{elastic-difference-farfield-effect-1**1-1-s}
{4\pi c^2_s} \left[ U^\infty_s(\hat{x},\theta)-U_{\bold{C}_{\textcolor{black}{\epsilon}},s}^\infty(\hat{x},\theta)\right]\cdot \hat{x}^{\perp}
      &=&O\left({\textcolor{black}{\epsilon}}^{\frac{1}{3}} +{\textcolor{black}{\epsilon}}^{\frac{3}{2}-3t}\right),\;~~ \frac{1}{3} \leq t < \frac{1}{2}.
  \end{eqnarray}
 \subsection{End of the proof of Theorem \ref{equivalent-medimu}}
  Combining the estimates (\ref{elastic-difference-farfield-effect-1**1-1}) and (\ref{appro-0-a}), we deduce that
\begin{equation}\label{final}
{4\pi c^2_p} \left[U_p^\infty(\hat{x}, \theta)-U_{0, p}^\infty(\hat{x}, \theta) \right ]\cdot \hat{x}=O({\textcolor{black}{\epsilon}}^{\min \{\gamma,\; \frac{1}{3},\; \frac{3}{2}-3t\}}),\; {\textcolor{black}{\epsilon}}<<1, \;~~ \frac{1}{3} \leq t < \frac{1}{2}
 \end{equation}
uniformly in terms of  $\hat x , \theta\; \in \mathbb{S}^{2}$.

Similarly, by combining the estimates (\ref{elastic-difference-farfield-effect-1**1-1-s}) and (\ref{appro-0-a}), we deduce that
\begin{equation}\label{final-s}
{4\pi c^2_s} \left[U_s^\infty(\hat{x}, \theta)-U_{0, s}^\infty(\hat{x}, \theta) \right ]\cdot \hat{x}^\perp=O({\textcolor{black}{\epsilon}}^{\min \{\gamma,\; \frac{1}{3},\; \frac{3}{2}-3t\}}),\; {\textcolor{black}{\epsilon}}<<1, \;~~ \frac{1}{3} \leq t < \frac{1}{2}
 \end{equation}
uniformly in terms of  $\hat x , \theta\; \in \mathbb{S}^{2}$.
 
 
\section{Justification of Corollary \ref{equivalent-medimu-rho}} 
 
For an obstacle $D_{\epsilon}$ of radius $\epsilon$, $\mathcal{S}(\phi)(s):=\int_{\partial D_{\epsilon}}\Gamma^\omega(s, t)\phi(t)dt$ and $\mathcal{D}(\phi)(s):=\int_{\partial D_{\epsilon}}\frac{\partial \Gamma^{\omega}(s, t)}{\partial \nu(t)}\phi(t)dt$. Similarly, we set $\mathcal{S}_G(\phi)(s):=\int_{\partial D_{\epsilon}}G_{\rho}(s, t)\phi(t)dt$ and $\mathcal{D}_G(\phi)(s):=\int_{\partial D_{\epsilon}}\frac{\partial G_{\rho}(s, t)}{\partial \nu(t)}\phi(t)dt$.
We see that $W_{\kappa}(x, z):=G_{\rho}(x, z)-\Gamma^{\omega}(x, z)$ satisfies
\begin{equation}
(\Delta +\omega^2)W_{\kappa}=\omega^2(1-\rho)\Gamma^{\omega},\; \mbox{ in } \mathbb{R}^3
\end{equation}
with the Kupradze radiation conditions.  Since $\Gamma^{\omega}(\cdot, z)$, $z \in \mathbb{R}^3$ is bounded in $L^p(\Omega)$, for $p<3$, by interior estimates, we deduce that $W(\cdot, z)$, $z \in \mathbb{R}^3$ is bounded in $W^{2, p}(\Omega)$, for $p<3$, and hence, in particular, the normal traces are bounded in $L^2(\partial D_{\epsilon})$. Then we can show that the norms of the operators
\begin{equation}\label{single-layer-rho}
\mathcal{S}_G- \mathcal{S}: (L^2(\partial D_{\epsilon}))^3 \rightarrow (H^1(\partial D_{\epsilon}))^3
\end{equation}   
and 
\begin{equation}\label{double-layer-rho}
\mathcal{D}_G- \mathcal{D}: (L^2(\partial D_{\epsilon}))^3 \rightarrow (L^2(\partial D_{\epsilon}))^3
\end{equation}
are of the order $O(\epsilon)$ at least. 
\begin{enumerate}
 
\item Using these properties and arguing as in \cite{C-S:2015}, we derive the asymptotic expansions (\ref{x oustdie1 D_m farmainp-near-rho})-(\ref{x oustdie1 D_m farmains-near-rho}). Indeed, apart from the computations done in \cite{C-S:2015}, the main arguments needed to extend those results to the case of variable density is the Fredholm alternative for the corresponding integral operators and the application of the Neumann series expansions. After splitting $G_{\rho}$ as $G_{\rho}=\Gamma^{\omega}+(G_{\rho}-\Gamma^{\omega})$, these two arguments are applicable as soon as we have (\ref{single-layer-rho})-(\ref{double-layer-rho}). 

\item The justification of the invertibility of the algebraic system (\ref{fracqcfracmain}) depends only on (1) the distribution of small bodies and (2) the background medium through the singularities of the fundamental solution (of the form $\vert \Gamma^{\omega}(s, t)\vert \leq c \vert s-t \vert^{-1}$). However, this type of singularity is true for general background elastic media \footnote{Of course, it can be justify using the decomposition $G_{\rho}=\Gamma^{\omega}+(G_{\rho}-\Gamma^{\omega})$ with the singularity of $\Gamma^{\omega}$ and the smoothness of $G_{\rho}-\Gamma^{\omega}$}. Then the same arguments can be used to justify the invertibiliy of the algebraic system (\ref{fracqcfracmain-rho}). Using the above mentioned decomposition of the Green's function $G_{\rho}$, the properties of the Lippmann Schwinger integral equation are also valid replacing $\Gamma^{\omega}$ by $G_{\rho}$, and hence the results in section \ref{subsection-piecewise-constant} are valid. Finally, and again using the decomposition of $G_{\rho}$, the computations in section \ref{the approximation-by-AS} can be carried out using $G_{\rho}$.       
 
\end{enumerate} 
 
\section{The elastic capacitance}\label{structure-capacitance}

We start with the following lemma on the symmetry structure of the elastic capacitance.
\begin{lemma}\label{Lemma-Cap-adj}
Let $C:=(\int_{\partial D}\sigma_{i, j}(t)dt)^3_{i, j=1}$ be the elastic capacitance of a bounded and Lischitz regular set $D$ and $C^*$ be its adjoint. Then
\begin{equation}\label{symmetry-capacitance}
C=C^*.
\end{equation}
\end{lemma}
\begin{proofe}{\it{Proof of Lemma \ref{Lemma-Cap-adj}.}}~ 
We know that the matrix $\sigma:=(\sigma_{i, j})^3_{i, j=1}$ solves the invertible integral equation $
\int_{\partial D}\Gamma_0(s, t)\sigma(t)dt=I_3,$ or precisely $\int_{\partial D}\Gamma_0(s, t)\sigma_i(t)dt=e_i$, where $\sigma_i:=(\sigma_{i, j})^3_{j=1}$ and $e_i$ is the $i^{th}$ column of $I_3$.
Let $a$ be any constant vector in $\mathbb{R}^3$, then the vector $\sigma\; a$ satisfies 
$$
\int_{\partial D}\Gamma_0(s, t)\left (\sigma(t)\; a\right)\; dt=a.
$$
We set $\varphi^a:=\int_{\partial D}\Gamma_0(s, t)\left (\sigma(t)\; a \right )dt$. Then $\varphi^a$ satisfies the problem $\Delta^e \varphi^a=0$, in $D$ and $\varphi^a=a$ on $\partial D$. In addition, we have the jump relation $\partial_{\nu +}\varphi^a-\partial_{\nu -}\varphi^a=\sigma(t)\; a$ on $\partial D$ where $\partial_\nu u:=\lambda(\nabla \cdot u)\nu +\mu(\nabla u +\nabla u^\top)\nu$ is the elastic conormal derivative. Hence
$$
\int_{\partial D}  \partial_{\nu +}\varphi^a-\partial_{\nu -}\varphi^a\; dt=\int_{\partial D}\sigma(t)\; a dt =C\;a.
$$
Now, let $a$ and $b$ be arbitrary constant vectors in $\mathbb{R}^3$. To both $a$ and $b$, we correspond $\varphi^a$ and $\varphi^b$ as above. Using the Green formulas inside and outside of $D$, we deduce that 
$$
\left(C\;a,b \right)=\int_{\partial D} \left (\partial_{\nu +}\varphi^a-\partial_{\nu -}\varphi^a(t)\right )\; \cdot \varphi^b(t)\; dt=\int_{\partial D}  \left( \partial_{\nu +}\varphi^b-\partial_{\nu -}\varphi^b \right )\cdot \varphi^a dt=\left(C\;b,a \right)=(a, C\;b)
$$
recalling that every quantity here is real valued.
\end{proofe}

The next lemma describes the elastic capacitance of a given bounded and Lipschitz regular domain with the one of its image by a unitary transform. 

\begin{lemma}\label{capacitance-rotation}
Let $\mathcal{R}=(r_{lm})$ be a unitary transform in $\mathbb{R}^d$, $D$ be bounded Lipschitz domain in $\mathbb{R}^d$, $d=2,3$ and $\tilde{D}=\mathcal{R}(D)$. 
Let $C$ and $\tilde{C}$ be the corresponding elastic capacitance matrices due to 
the density matrices $\sigma$ and $\tilde{\sigma}$, as defined in \eqref{barqcimsurfacefrm1main}, respectively.
Then we have $\tilde{C}=\mathcal{R}\;C\;\mathcal{R}^{-1}$.
\end{lemma}
\begin{proofe}{\it{Proof of Lemma \ref{capacitance-rotation}.}}~ First recall the relation $\Gamma^{0}\circ\mathcal{R}(\xi,\eta)= \mathcal{R}\Gamma^{0}(\xi,\eta)\mathcal{R}^{-1}$, see \cite[Lemma 6.11]{AH-KH:SV18462004}. From \eqref{barqcimsurfacefrm1main}, we have that
\begin{eqnarray}\label{barqcimsurfacefrm1main-1}
&\int_{\partial \tilde{D}}\Gamma^{0}(\tilde{\xi},\tilde{\eta})\tilde{\sigma} (\tilde{\eta})d\tilde{\eta}
&\,=\,\rm \textbf{I},~ \tilde{\eta}\in \partial \tilde{D}\nonumber\\
\Rightarrow&\int_{\partial D}\left(\Gamma^{0}\circ\mathcal{R}\right)(\xi,\eta)\left(\tilde{\sigma}\circ\mathcal{R}\right) (\eta)d\eta
&\,=\,\rm \textbf{I},~ {\eta}\in \partial {D}\nonumber\\
\Rightarrow&\int_{\partial D}\mathcal{R}\Gamma^{0}(\xi,\eta)\mathcal{R}^{-1}\tilde{\sigma}\circ\mathcal{R} (\eta)d\eta
&\,=\,\rm \textbf{I},~ {\eta}\in \partial {D}\nonumber\\
\Rightarrow&\int_{\partial D}\Gamma^{0}(\xi,\eta)\mathcal{R}^{-1}\left(\tilde{\sigma}\circ\mathcal{R}\right) (\eta)d\eta
&\,=\,\mathcal{R}^{-1},~ {\eta}\in \partial {D}\nonumber\\
\Rightarrow&\int_{\partial D}\Gamma^{0}(\xi,\eta)\mathcal{R}^{-1}\left(\tilde{\sigma}\circ\mathcal{R}\right) (\eta)\mathcal{R}d\eta
&\,=\,\rm \textbf{I},~ {\eta}\in \partial {D}.
\end{eqnarray}
Now from the uniqueness of solutions of \eqref{barqcimsurfacefrm1main}, we deduce that
$\mathcal{R}^{-1}\left(\tilde{\sigma}\circ\mathcal{R}\right)(\cdot) \mathcal{R}\,=\,\sigma(\cdot)$ and then
 \begin{eqnarray}\label{density-rotation}
\left(\tilde{\sigma}\circ\mathcal{R}\right)(\cdot)\,=\, \mathcal{R}\sigma(\cdot)  \mathcal{R}^{-1}.
\end{eqnarray}
From the definition of the capacitance, see (\ref{barqcimsurfacefrm1main}), and \eqref{density-rotation}, we have
\begin{eqnarray}\label{cap-rotation}
\tilde{C}
=\int_{\partial \tilde{D}}\tilde{\sigma}(\tilde{\eta})d\tilde{\eta}
=\int_{\partial {D}}\left(\tilde{\sigma}\circ\mathcal{R}\right) (\eta)d{\eta}
=\int_{\partial {D}}\mathcal{R}\sigma(\eta)  \mathcal{R}^{-1} d{\eta}
=\mathcal{R}{C}\mathcal{R}^{-1}.
\end{eqnarray}
\end{proofe}
\begin{proposition}\label{capacitance-rotation-2Drltn}
Let $\mathcal{R}=(r_{lm})$ be a unitary transform in $\mathbb{R}^d$, $d=2, 3$, and 
let $D$ be a bounded Lipschitz domain in $\mathbb{R}^d$, $d=2,3$, and $\tilde{D}=\mathcal{R}(D)$. 
Let $C$ and $\tilde{C}$ be the corresponding elastic capacitance matrices due to 
the density matrices $\sigma$ and $\tilde{\sigma}$, as defined in \eqref{barqcimsurfacefrm1main}, respectively.
\begin{enumerate}
\item $2D$-case. If the shape of $D$ is rotationally invariant for any rotation by one angle $\theta \neq 0, \pi$,
then $C$ is a scalar multiplied by the identity matrix.

\item $3D$-case. If the shape of $D$ is rotationally invariant for any two of the rotations around the $x, y$ or $z$ axis by  one angle $\theta \neq 0, \pi$ and $\alpha \neq 0, \pi$ respectively,
then $C$ is a scalar multiplied by the identity matrix.
\end{enumerate}
\end{proposition}

\begin{proofe}{\it{Proof of Proposition \ref{capacitance-rotation-2Drltn}.}}~ 
 In $2\textbf{D}$ case, the rotation matrix by an angle $\theta$ is given by
 \begin{eqnarray}\label{Rotationmatrix-2D}
  \mathcal{R}&=&\left(\begin{array}{ccc}
                       \cos\theta&-\sin\theta\\
                       \sin\theta&\cos\theta
                      \end{array}\right).
 \end{eqnarray}
As the shape is invariant by this rotation then $\tilde{C}=C$.
Since $\mathcal{R}$ is unitary then $\mathcal{R}^{-1}=\mathcal{R}^{\top}$ and then \eqref{cap-rotation} implies
\begin{eqnarray}\label{cap-rotation-1}
C&&\,=\tilde{C}=\mathcal{R}{C}\mathcal{R}^{\top}\nonumber\\
=&&\left(\begin{array}{ccc}
                       \cos\theta&-\sin\theta\\
                       \sin\theta&\cos\theta
                      \end{array}\right)
                      \left(\begin{array}{ccc}
                       C_{11}&C_{12}\\
                       C_{21}&C_{22}
                      \end{array}\right)
                      \left(\begin{array}{ccc}
                       \cos\theta&\sin\theta\\
                       -\sin\theta&\cos\theta
                      \end{array}\right)\nonumber\\
=&&\left(\begin{array}{ccc}
                       \cos\theta&-\sin\theta\\
                       \sin\theta&\cos\theta
                      \end{array}\right)
                      \left(\begin{array}{ccc}
                       C_{11}\cos\theta-C_{12}\sin\theta&C_{11}\sin\theta+C_{12}\cos\theta\\
                       C_{21}\cos\theta-C_{22}\sin\theta&C_{21}\sin\theta+C_{22}\cos\theta
                     \end{array}\right)\\
=&&\left(\begin{array}{ccc}
                       C_{11}\cos^2\theta+C_{22}\sin^2\theta-
                       (C_{12}+C_{21})\sin\theta\cos\theta
                       &C_{12}\cos^2\theta-C_{21}\sin^2\theta+(C_{11}-C_{22})\sin\theta\cos\theta\\
                       C_{21}\cos^2\theta-C_{12}\sin^2\theta+(C_{11}-C_{22})\sin\theta\cos\theta
                       &C_{11}\sin^2\theta+C_{22}\cos^2\theta+
                       (C_{12}+C_{21})\sin\theta\cos\theta
                     \end{array}\right). \nonumber                   
\end{eqnarray}

We deduce from (\ref{cap-rotation-1}) and the symmetry of matrix $C$ the following relations:
\begin{equation}\label{relations-2D-1}
C_{11}=C_{11}\cos^2\theta+C_{22}\sin^2\theta-
                       2C_{12}\sin\theta\cos\theta
\end{equation}
\begin{equation}\label{relations-2D-2}
C_{22}=C_{11}\sin^2\theta+C_{22}\cos^2\theta+
                       2C_{12}\sin\theta\cos\theta
\end{equation}
\begin{equation}\label{relations-2D-3}                                              
C_{12}=C_{12}\cos^2\theta-C_{12}\sin^2\theta+(C_{11}-C_{22})\sin\theta\cos\theta.
\end{equation}
We rewrite (\ref{relations-2D-1}) and (\ref{relations-2D-3}) respectively as
\begin{equation}\label{relations-2D-1-1}
(C_{11}-C_{22})\sin^2\theta+
                       2C_{12}\sin\theta\cos\theta=0
\end{equation}
\begin{equation}\label{relations-2D-3-1}
(C_{11}-C_{22})\cos\theta \sin \theta-
                       2C_{12}\sin^2\theta=0
\end{equation}
Taking $\sin \theta$ as a multiplicative factor in (\ref{relations-2D-1-1}) and (\ref{relations-2D-3-1}) we see that if $\theta \neq 0, \pi$, i.e. $\sin \theta \neq 0$, then we have
\begin{equation}\label{2D-conlcusion}
C_{11}=C_{22} \mbox{ and } C_{12}=C_{21}=0.
\end{equation}
\bigskip

Let us now consider the 3D case. 
First, let us assume that the shape is invariant under the rotation about the $x-axis$ and an angle $\theta \neq 0, \pi$. This rotation matrix is given by
 \begin{eqnarray}\label{Rotationmatrix-3D-x}
  \mathcal{R}&=&\left(\begin{array}{ccc}
                     1 &       0     &  0        \\
                     0 &  \cos\theta & -\sin\theta\\
                     0 & \sin\theta  & \cos\theta
                      \end{array}\right)
 \end{eqnarray}
Since $\mathcal{R}$ is unitary, $\mathcal{R}^{-1}=\mathcal{R}^{\top}$, then \eqref{cap-rotation} gives us;
\begin{eqnarray}\label{cap-rotation-1-3D-x}
C=\tilde{C}&=&\mathcal{R}{C}\mathcal{R}^{\top}\nonumber\\
&=&\left(\begin{array}{ccc}
                      1 &       0     &  0        \\
                     0 &  \cos\theta & -\sin\theta\\
                     0 & \sin\theta  & \cos\theta
                      \end{array}\right)
                      \left(\begin{array}{ccc}
                       C_{11}&C_{12}&C_{13}\\
                       C_{21}&C_{22}&C_{23}\\
                       C_{31}&C_{32}&C_{33}
                      \end{array}\right)
                      \left(\begin{array}{ccc}
                       1 &       0     &  0        \\
                     0 &  \cos\theta & \sin\theta\\
                     0 & -\sin\theta  & \cos\theta
                      \end{array}\right)\nonumber\\
&=&\left(\begin{array}{ccc}
                       1 &       0     &  0        \\
                     0 &  \cos\theta & -\sin\theta\\
                     0 & \sin\theta  & \cos\theta
                      \end{array}\right)
                      \left(\begin{array}{ccc}
                      C_{11} & C_{12}\cos\theta-C_{13}\sin\theta & C_{12}\sin\theta+C_{13}\cos\theta\\
                      C_{21} & C_{22}\cos\theta-C_{23}\sin\theta & C_{22}\sin\theta+C_{23}\cos\theta\\
                      C_{31} & C_{32}\cos\theta-C_{33}\sin\theta & C_{32}\sin\theta+C_{33}\cos\theta
                     \end{array}\right)\\
&=&\left(\begin{array}{ccc}
                         T_{11} & T_{12} & T_{13}
                         \\
                         T_{21} & T_{22} & T_{23}
                         \\
                       T_{31} & T_{32}   & T_{33}
                     \end{array}\right),
                     \hspace{1.5cm} \nonumber                   
\end{eqnarray}
 \begin{eqnarray}\label{defining-Tsubsc}
 \mbox{where\qquad }&
\left\{\begin{array}{ccc}
                        T_{11} &:=& C_{11};\hspace{6.5cm}\\
                        T_{12} &:=& C_{12}\cos\theta-C_{13}\sin\theta;\hspace{4cm} \\ 
                        T_{13} &:=& C_{12}\sin\theta+C_{13}\cos\theta;\hspace{4cm}
                         \\
                         T_{21} &:=& C_{21}\cos\theta-C_{31}\sin\theta;\hspace{4cm}\\
                         T_{22} &:=& C_{22}\cos^2\theta+C_{33}\sin^2\theta-(C_{22}+C_{33})\sin\theta\cos\theta;\\
                         T_{23} &:=& C_{23}\cos^2\theta-C_{32}\sin^2\theta+(C_{22}-C_{33})\sin\theta\cos\theta;
                       \\
                        T_{31} &:=& C_{21}\sin\theta+C_{31}\cos\theta;\hspace{4cm}\\
                        T_{32} &:=& C_{32}\cos^2\theta-C_{23}\sin^2\theta+(C_{22}-C_{33})\sin\theta\cos\theta;\\
                        T_{33} &:=& C_{22}\sin^2\theta+C_{33}\cos^2\theta+(C_{23}+C_{32})\sin\theta\cos\theta;
 \end{array}\right.& \hspace{4.85cm}
 \end{eqnarray}

We observe the equality of the $2 \times 2$ matrices:
 \begin{align} \label{cap-rotation-1-3D-x-1}
    &\left(\begin{array}{ccc}                      
           C_{22}&C_{23}\\
           C_{32}&C_{33}
          \end{array}\right) = \left(\begin{array}{ccc}                      
           T_{22}&T_{23}\\
           T_{32}&T_{33}
          \end{array}\right)  \\                  
 = &\left(\begin{array}{ccc}                     
       C_{22}\cos^2\theta+C_{33}\sin^2\theta-(C_{22}+C_{33})\sin\theta\cos\theta
     & C_{23}\cos^2\theta-C_{32}\sin^2\theta+(C_{22}-C_{33})\sin\theta\cos\theta
       \\
       C_{32}\cos^2\theta-C_{23}\sin^2\theta+(C_{22}-C_{33})\sin\theta\cos\theta
     & C_{22}\sin^2\theta+C_{33}\cos^2\theta+(C_{23}+C_{32})\sin\theta\cos\theta
 \end{array}\right).\hspace{1.5cm} \nonumber                   
\end{align}

These are similar to the matrices we obtained in the $2D$ case. Hence we deduce, as in the $2D$ case, that
\begin{equation}\label{C2-s}
C_{22}=C_{33} \mbox{ and } C_{23}=C_{32}=0.
\end{equation} 

To show that $C$ is scalar multiplied by the identity matrix we need to prove that $C_{11}=C_{22}$, for instance, and $C_{13}=C_{31}=0$.  
For this purpose, we use another rotation. Taking the rotation around the $z$-axis \footnote{ We can use also the rotation around the $y$ axis.} by one angle $\alpha \neq 0,\pi$ and proceeding as we did for the rotation about the $x$-axis,
we show that
\begin{equation}\label{C1-s}
C_{11}=C_{33} \mbox{ and } C_{13}=C_{31}=0.
\end{equation} 
 \end{proofe}
 \bigskip

 From the above analysis, we have the following remark:
 \begin{remark}\label{cap-spherical}
\begin{enumerate}
\item For the spherical shapes, in particular, the capacitance is a scalar multiplied by the identity matrix. 

\item Ellipsoidal shapes are invariant only under rotations with angle $\pi$ (or trivially $0$). For these shapes, the capacitance might not be a scalar multiplied by the identity matrix but a diagonal matrix instead. To justify this property, the arguments in \cite{AH-KH:SV18462004} can be useful. 
\end{enumerate}
  \end{remark} 
 



\bibliographystyle[numbib]{imamat}
%



\begin{thebibliography}{}

\bibitem[Ahmad et~al., 2015]{DPC-SM13-3}
Ahmad, B., Challa, D.~P., Kirane, M. {\&} Sini, M. (2015)  The equivalent
  refraction index for the acoustic scattering by many small obstacles: with
  error estimates. {\em J. Math. Anal. Appl.}, \textbf{424}(1), 563--583.

\bibitem[Alves \& Kress, 2002]{A-K:IMA2002}
Alves, C. J.~S. {\&} Kress, R. (2002)  On the far-field operator in elastic
  obstacle scattering. {\em IMA J. Appl. Math.}, \textbf{67}(1), 1--21.

\bibitem[Ammari \& Kang, 2004]{AH-KH:SV18462004}
Ammari, H. {\&} Kang, H. (2004) {\em Reconstruction of small inhomogeneities
  from boundary measurements}, volume 1846 of {\em Lecture Notes in
  Mathematics}.
Springer-Verlag, Berlin.

\bibitem[Ammari et~al., 2007]{A-K-L-LameCPDE2007}
Ammari, H., Kang, H. {\&} Lee, H. (2007)  Asymptotic expansions for eigenvalues
  of the {L}am\'e system in the presence of small inclusions. {\em Comm.
  Partial Differential Equations}, \textbf{32}(10-12), 1715--1736.

\bibitem[Bensoussan et~al., 1978]{B-L-P:1978}
Bensoussan, A., Lions, J.-L. {\&} Papanicolaou, G. (1978) {\em Asymptotic
  analysis for periodic structures}, volume~5 of {\em Studies in Mathematics
  and its Applications}.
North-Holland Publishing Co., Amsterdam.

\bibitem[Challa \& Sini, 2014]{DPC-SM13}
Challa, D.~P. {\&} Sini, M. (2014)  On the justification of the {F}oldy-{L}ax
  approximation for the acoustic scattering by small rigid bodies of arbitrary
  shapes. {\em Multiscale Model. Simul.}, \textbf{12}(1), 55--108.

\bibitem[Challa \& Sini, 2015]{C-S:2015}
Challa, D.~P. {\&} Sini, M. (2015)  The Foldy-Lax approximation of the
  scattered waves by many small bodies for the Lamé system. {\em Mathematische
  Nachrichten}, \textbf{288}(16), 1834--1872.

\bibitem[Cioranescu \& Murat, 1982]{DC-FM:Book:BB1979}
Cioranescu, D. {\&} Murat, F. (1982)  Un terme \'etrange venu d'ailleurs. In
  {\em Nonlinear partial differential equations and their applications.
  {C}oll\`ege de {F}rance {S}eminar, {V}ol. {II} ({P}aris, 1979/1980)},
  volume~60 of {\em Res. Notes in Math.}, pages 98--138, 389--390. Pitman,
  Boston, Mass.-London.

\bibitem[Cioranescu \& Murat, 1997]{DC-FM:Book:BB1997}
Cioranescu, D. {\&} Murat, F. (1997)  A strange term coming from nowhere [
  {MR}0652509 (84e:35039a); {MR}0670272 (84e:35039b)]. In {\em Topics in the
  mathematical modelling of composite materials}, volume~31 of {\em Progr.
  Nonlinear Differential Equations Appl.}, pages 45--93. Birkh\"auser Boston,
  Boston, MA.

\bibitem[Colton \& Kress, 1998]{C-K:1998}
Colton, D. {\&} Kress, R. (1998) {\em Inverse acoustic and electromagnetic
  scattering theory}, volume~93 of {\em Applied Mathematical Sciences}.
Springer-Verlag, Berlin, second edition.

\bibitem[Colton \& Kress, 1983]{C-K:1983}
Colton, D.~L. {\&} Kress, R. (1983) {\em Integral equation methods in
  scattering theory}.
Pure and Applied Mathematics (New York). John Wiley \& Sons Inc., New York.
A Wiley-Interscience Publication.

\bibitem[Hu \& Liu, 2015]{HG-LH:JMPANS2015}
Hu, G. {\&} Liu, H. (2015)  Nearly cloaking the elastic wave fields. {\em J.
  Math. Pures Appl. (9)}, \textbf{104}(6), 1045--1074.

\bibitem[Jikov et~al., 1994]{J-K-O:1994}
Jikov, V.~V., Kozlov, S.~M. {\&} Ole{\u\i}nik, O.~A. (1994) {\em Homogenization
  of differential operators and integral functionals}.
Springer-Verlag, Berlin.

\bibitem[Kupradze, 1965]{Kupradze:1965}
Kupradze, V.~D. (1965) {\em Potential methods in the theory of elasticity}.
Translated from the Russian by H. Gutfreund. Translation edited by I. Meroz.
  Israel Program for Scientific Translations, Jerusalem.

\bibitem[Kupradze et~al., 1979]{K-G-B-B:1979}
Kupradze, V.~D., Gegelia, T.~G., Bashele{\u\i}shvili, M.~O. {\&} Burchuladze,
  T.~V. (1979) {\em Three-dimensional problems of the mathematical theory of
  elasticity and thermoelasticity}, volume~25 of {\em North-Holland Series in
  Applied Mathematics and Mechanics}.
North-Holland Publishing Co., Amsterdam, russian edition.
Edited by V. D. Kupradze.

\bibitem[Marchenko \& Khruslov, 2006]{M-K:2006}
Marchenko, V.~A. {\&} Khruslov, E.~Y. (2006) {\em Homogenization of partial
  differential equations}, volume~46 of {\em Progress in Mathematical Physics}.
Birkh{\"a}user Boston Inc., Boston, MA.

\bibitem[McLean, 2000]{Mclean:2000}
McLean, W. (2000) {\em Strongly elliptic systems and boundary integral
  equations}.
Cambridge University Press, Cambridge.

\bibitem[Namias, 1986]{NV:Stirling:AMM:1986}
Namias, V. (1986)  A simple derivation of {S}tirling's asymptotic series. {\em
  Amer. Math. Monthly}, \textbf{93}(1), 25--29.


\bibitem[Ramm, 2007]{RAMM:2007}
Ramm, A.~G. (2007)  Many-body wave scattering by small bodies and applications.
  {\em J. Math. Phys.}, \textbf{48}(10), 103511, 29.

\bibitem[Ramm, 2011]{RAMM:2011}
Ramm, A.~G. (2011)  Wave scattering by small bodies and creating materials with
  a desired refraction coefficient. {\em Afr. Mat.}, \textbf{22}(1), 33--55.

\end{thebibliography}

\def\cprime{$'$}

\end{document}